%% file: Paper_v2_double-column.tex
\documentclass[10pt,twocolumn,twoside]{IEEEtran}

\input{commands}
\usepackage[dvips]{graphicx}
\usepackage{amsmath}
\usepackage{subfigure}
\usepackage{cite}
\usepackage{url}
\usepackage{color}
\usepackage{xcolor}
\usepackage{amssymb}
\usepackage{psfrag}
\usepackage{caption}
\usepackage{booktabs}
\usepackage{array}
\usepackage{epstopdf}
\usepackage{epsfig}
\usepackage{tikz}
\usepackage{lipsum}
\usetikzlibrary{arrows,decorations.pathmorphing,backgrounds,fit,petri,positioning,matrix}
\usetikzlibrary{shapes,arrows}
\usepackage{multirow}
\renewcommand{\thetable}{\arabic{table}}
\usetikzlibrary{patterns}
\newlength{\hatchspread}
\newlength{\hatchthickness}
\tikzset{hatchspread/.code={\setlength{\hatchspread}{#1}},
         hatchthickness/.code={\setlength{\hatchthickness}{#1}}}
\tikzset{hatchspread=3pt,
         hatchthickness=0.4pt}
\pgfdeclarepatternformonly[\hatchspread,\hatchthickness]
   {custom north west lines}
   {\pgfqpoint{-2\hatchthickness}{-2\hatchthickness}}
   {\pgfqpoint{\dimexpr\hatchspread+2\hatchthickness}{\dimexpr\hatchspread+2\hatchthickness}}
   {\pgfpoint{\hatchspread}{\hatchspread}}
   {
    \pgfsetlinewidth{\hatchthickness*rand}
    \pgfpathmoveto{\pgfpoint{rand*0.2pt}{\hatchspread}}
   \pgfpathcurveto
{\pgfqpoint{\dimexpr\hatchspread+6pt}{0.2pt}}{\pgfpoint{\hatchspread+4pt}{rand*3pt}}{\pgfqpoint{\dimexpr\hatchspread+0.1pt}{0.15pt}}
    \pgfsetstrokeopacity{0.175}
    \pgfsetstrokecolor{blue}
    \pgfusepath{stroke}
   }

\usepackage{tkz-euclide}
\usetkzobj{all}

\newcommand{\addw}[1]{{{\color{white!100!black}#1}}}
\newcommand{\addb}[1]{{{\color{blue!60!white}#1}}}
\newcommand{\goodgap}{\hspace{0.01\linewidth}} 

\def\WiT#1{\color{white}{#1}}

\begin{document} 
\title{A Semi Distributed Approach for Min-Max Fair Car-Parking Slot Assignment Problem}
\author{E. Alfonsetti,
        P. C. Weeraddana,~\IEEEmembership{Member,~IEEE,} and
        C.~Fischione,~\IEEEmembership{Member,~IEEE}
\thanks{Manuscript received ... This research was supported by EU projects Hycon2, Hydrobionets, and VR project In network Optimization.}
\thanks{E. Alfonsetti with TerraSwarm Lab, Electrical Engineering Department, UC Berkeley, California, USA (e-mail: e.alfonsetti@berkeley.edu).}
\thanks{P. C. Weeraddana and C. Fischione with Electrical Engineering, KTH Royal Institute of Technology, Stockholm, Sweden (e-mail: chatw@kth.se, carlofi@kth.se).}
}

\maketitle

\begin{abstract}
Designing efficient car parking mechanisms that can be potentially integrated
into future intelligent transportation systems is of crucial importance. Usually, the related design problems are combinatorial and the worst-case complexity of optimal solution approaches grows exponentially with the problem sizes. Therefore, such optimal approaches are not scalable and practically undesirable. As a result, almost all existing methods for parking slot assignment are simple and greedy approaches, where each car is assigned a free parking slot, which is closer to its destination. Moreover, no emphasis is placed to optimize the social benefit of the users during the parking slot assignment. In this paper, the fairness as a metric for modeling the aggregate social benefit of the users is considered and a distributed algorithm based on Lagrange duality theory is developed. The proposed algorithm is gracefully scalable compared to the optimal methods. In addition, it is shown that the proposed car parking mechanism preserves privacy in the sense that any car involved in the algorithm will not be able to discover the destination of any other car during the algorithm iterations. Numerical results illustrate the performance of the proposed algorithm compared to the optimal assignment and a greedy method. They show that our algorithm yields a good tradeoff between the implementation-level simplicity and the~performance. Even though the main emphasis in this paper resides in the car parking slot assignment problem, our formulation and the algorithms, in general, can also be applied or adopted in fair agent-target assignment problems in other application domains.
\end{abstract}

\vspace{-0mm}
\begin{keywords}\vspace{-0mm}
Intelligent transportation systems, optimization methods, algorithms, privacy
\end{keywords}


\section{Introduction} \label{sec:intro}

The car is certainly one of the most used means of transport, but its introduction, despite having brought comfort and simplification of life, has generated well-known problems of increasing traffic, and thus clogging roads in the city centers. Especially in large cities, these problems are pronounced by hundreds or even thousands of drivers who are looking for parking slots during their daily activities. In~\cite{Rodrigue-Comtois-Slack-2013} it is claimed that seeking for parking slots (\emph{cruising}) can account for more than 10\% of the local circulations in central areas of large cities. In \cite{Shoup-2006}, it is reported that cruising for open parking spaces accounts for 30\% of the traffic, causing undesired congestion in big cities. In addition, cruising creates additional delays and drivers can even spend up to 10-20 minutes before they could find a proper parking slot. According to a recent British study, it is estimated that a person who owns a car in big cities can take an average of 6 to 20 minutes to search for an empty parking slot, which accounts for monetary losses (e.g., deterioration, unnecessary fuel wastage), as well as for nonmonetary expenses (e.g., frustration, psychophysical stresses).

In the context of future intelligent transportation systems (ITS), there are many relevant research activities, which design various traffic congestion control mechanisms, see \cite{Chen-Cheng-2010,Jianming-Qiang-Qi-Jiajie-Yi-2012,Lee-etal-2011,Dimitrakopoulos-etal-2012,Chiew-Shaowen-2009} and references therein. However, \cite{Shoup-2006} suggests that designing efficient car parking mechanisms, which are instrumental in directly reducing the cruising traffic is just as important as other related methodologies to minimize undesired traffic conditions, especially, in big cities.



Several research attempts have been made in the field of ITS, which support drivers to locate a free parking slots, see \cite{Ergen-et-al-2003 ,Grossglauser-Piorkowski-2006,Piorkowski-etal-MobiHoc-2006,Piorkowski-etal-TechDemo-2006,Tang-Zheng-Cao-2006,Bi-etal-2006,Chinrungrueng-Sunantachaikul-etal-2007,Rongxing-etal-2009,Rongxing-etal-2010,Souissi-etal-2011}. In general, these existing methods employ a central authority (CA) who is responsible for providing the underlying infrastructure. For example, thousands of sensors have to be deployed to detect the availability of free parking slots~\cite{Ergen-et-al-2003 ,Grossglauser-Piorkowski-2006,Piorkowski-etal-MobiHoc-2006,Piorkowski-etal-TechDemo-2006,Tang-Zheng-Cao-2006,Bi-etal-2006,Chinrungrueng-Sunantachaikul-etal-2007,Souissi-etal-2011}. The authors in \cite{Rongxing-etal-2009,Rongxing-etal-2010} have considered a stage of vehicular ad hoc networks, where advanced roadside units are widely deployed and every vehicle is equipped with sophisticated onboard units. In addition, efficient database management systems together with real-time communication protocols have to be implemented to meet the demands.

The proliferation of smartphones has encouraged a number of private, as well as state companies to deploy real-time car parking mechanisms in central areas of large cities, see ParkingCarma~\cite{ParkingCarma}, Steeteline Parker~\cite{StreetLine}, VehicleSense Street Parking Information Network (SPIN)~\cite{VehicleSense}, VehicleSense SmartLot~\cite{VehicleSense}, and SFMTA SFpark~\cite{SFMTA}. Like other research work proposed in\cite{Ergen-et-al-2003 ,Grossglauser-Piorkowski-2006,Piorkowski-etal-MobiHoc-2006,Piorkowski-etal-TechDemo-2006,Tang-Zheng-Cao-2006,Bi-etal-2006,Chinrungrueng-Sunantachaikul-etal-2007,Rongxing-etal-2009,Rongxing-etal-2010,Souissi-etal-2011}, these methods are also relying on a central authority for providing the underlying infrastructure, such as wireless sensors, database management systems, etc.

\begin{figure}[t]
\centering

\subfigure[]{
\begin{tikzpicture}[scale=2.0]
\tikzstyle{state}=[circle,thick,draw=black!75,fill=black!20,minimum size=6mm, text=black]
\tikzstyle{ParkSlot}=[rectangle,rotate=90,draw=black,fill=blue!60!white,thick,inner sep=0pt,minimum size=3.0mm]
\tikzstyle{DesNode}=[rectangle,rotate=90,draw=black!100!white,fill=black!60!white,thick,inner sep=0pt,minimum size=2.0mm]
\tikzstyle{MyMatrix1} = [matrix, draw=blue, thick, fill=blue!20, text centered,font=\footnotesize,minimum height=2em]
\tikzstyle{agent1}=[circle,draw=red,fill=white, inner sep=0pt,minimum size=1.0mm]
\tikzstyle{agent2}=[circle,draw=blue,fill=blue, inner sep=0pt,minimum size=1.0mm]
\tikzstyle{target1}=[rectangle,draw=black,fill=white, inner sep=0pt,minimum size=1mm]
\tikzstyle{target2}=[rectangle,draw=black,fill=black, inner sep=0pt,minimum size=1mm]

\draw[color=blue, thick, opacity=0.2] (0,0.5) circle (0.04\textwidth);
\draw[color=red, thick, opacity=0.2] (0, 0.5) circle (0.01\textwidth);

\draw[color=blue, thick, opacity=0.2] (0.75,0.5) circle (0.05\textwidth);
\draw[color=red, thick, opacity=0.2] (0.75, 0.5) circle (0.04\textwidth);

\node at (0,0.5) [target1] (Xpoint1) {\WiT{.}};
\node at (0.75,0.5) [target2] (Xpoint1) {.};

\node at (0.05,0.67) [agent1] (Xpoint1) {\WiT{.}};
\node at (0.18,1.19) [agent2] (Xpoint1) {.};

\end{tikzpicture}
\label{fig:MinmaxVsSum-b}}
\quad\quad\quad\quad
\subfigure[]{
\begin{tikzpicture}[scale=1.4]
\tikzstyle{state}=[circle,thick,draw=black!75,fill=black!20,minimum size=6mm, text=black]
\tikzstyle{ParkSlot}=[rectangle,rotate=90,draw=black,fill=blue!60!white,thick,inner sep=0pt,minimum size=3.0mm]
\tikzstyle{DesNode}=[rectangle,rotate=90,draw=black!100!white,fill=black!60!white,thick,inner sep=0pt,minimum size=2.0mm]
\tikzstyle{MyMatrix1} = [matrix, draw=blue, thick, fill=blue!20, text centered,font=\footnotesize,minimum height=2em]
\tikzstyle{agent1}=[circle,draw=red,fill=white, inner sep=0pt,minimum size=3.0mm]
\tikzstyle{cloud} = [draw=white, node distance=4em,font=\small]
\tikzstyle{agent2}=[circle,draw=blue,fill=blue, inner sep=0pt,minimum size=3.0mm]
\tikzstyle{target1}=[rectangle,draw=black,fill=white, inner sep=0pt,minimum size=3mm]
\tikzstyle{target2}=[rectangle,draw=black,fill=black, inner sep=0pt,minimum size=3mm]

\node at (1,1.35) [agent1] (a1) {} ;
\node at (1,0) [agent2] (a2) {} ;
\node at (2,1.35) [target1] (t1) {} ;
\node at (2,0) [target2] (t2) {} ;

\node [cloud, left of=a1] (agentt1) {{car or user $1$ }};
\node [cloud, left of=a2] (agentt1) {{car or user $2$}};
\node [cloud, right of=t1] (tt1) {{parking slot $1$}};
\node [cloud, right of=t2] (tt1) {{parking slot $2$}};

\draw (a1) edge [-latex, above, thick] node {$1$} (t1);
\draw (a1) edge [-latex, near start, above, thick] node {$4$} (t2);
\draw (a2) edge [-latex, near start, above, thick] node {$4$} (t1);
\draw (a2) edge [-latex, above, thick] node {$5$} (t2);

\end{tikzpicture}
\label{fig:MinmaxVsSum-a}}

\captionof{figure}{\label{fig:MinmaxVsSum}\small{Car-Parking slot assignment: (a) Physical locations of free parking slots (squares) and cars' destinations (circles); (b) Graph, where the weights denote the distances between the destinations and parking slots.}}
\vspace{-5mm}
\end{figure}
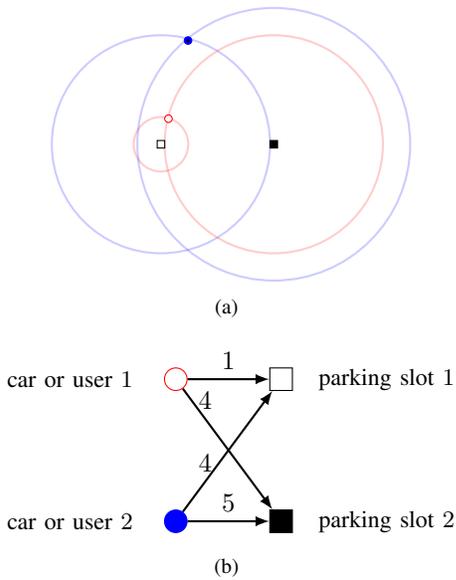
In almost all existing methods, no emphasis has been placed to optimize an aggregate social benefit of the users during the parking slot assignment. Arguably, the existing mechanisms can be interpreted as greedy methods, where each user selects the closest free parking slot to its destination. Such a greedy method can easily account for substantial imbalances among the distances between users' parking slots and their destinations, i.e.,~is not fair. In other words, some users can be assigned to parking slots that are very close to their destinations while others can be assigned to parking slots that are far from their destinations. For example, Figure~\ref{fig:MinmaxVsSum} shows a case where two cars are to be assigned to two parking slots. A greedy assignment approach yields the (car, parking slot) assignment $(1,1)$ and $(2,2)$ which accounts for a total cost of $6$ (i.e., $1+5$) and a cost imbalance of $5$. On the other hand, a fair assignment yields the assignment $(1,2)$ and $(2,1)$ which accounts for a cost imbalance of $1$.

Therefore, it is worth seeking efficient, as well as advanced algorithms that are capable of optimizing some aggregate social benefit for the users of the system. Because the involvement of a central authority is instrumental in coordinating the car parking mechanisms, optimization criteria can be integrated in to the parking assignment methods, where some aggregate social metric (utility) is considered during the assignment process. One such simple, yet appealing utility is users' fairness~\cite{Radunovic-Boudec-2007}.

Ensuring privacy of the associated algorithms, which can be in various contexts, is also important, see \cite{Rongxing-etal-2010}. Naturally, users would not like to publish information, such as their destinations to prevent a third party from predicting private traveling patterns. For example, government agencies can probe such information during investigations and business entities might be interested in exploiting such information to promote their products and services. Therefore, exposure of private information raises serious concerns of personal~privacy.

The main contributions of this paper are as follows:
\begin{enumerate}
\item We consider the \emph{min-max fairness} as a metric for modeling the aggregate social benefit of the users~\cite{Radunovic-Boudec-2007}. In particular, we consider the distance between parking slot and the destination that corresponds to every user. We refer to this distance, associated with any user, as the \emph{parking distance}. Then we design an algorithm to minimize the maximum parking distance among all the users. The proposed algorithm is based on duality theory~\cite[Section~5]{Boyd-Vandenberghe-04}. Our formulation and the corresponding algorithm can be applied directly or with minor modifications in fair agent-target assignment problems in other application domains as well, and therefore is not restricted to the car-parking assignment.

\item  We capitalize on dual decomposition techniques~\cite{Boyd-EE364b-PrimDualDecomp-07} and the subgradient methods~\cite{Boyd-EE364b-SubGradMethods-07} to accomplish distributed implementation (among users) of the proposed algorithm with a little coordination of the central authority. Therefore, the proposed algorithm has rich scalability properties, which is indeed favorable in practice.

\item The proposed car parking mechanism is privacy preserving in the sense that any
car involved in the algorithm will not be able to find out the destination of any other car during the algorithm iterations. This privacy is accomplished as a result of the inherent decomposition structure of the problem together with randomization of the step size of the subgradient method.

\item  A number of numerical examples are provided to evaluate the performance of the algorithm. In addition, the proposed algorithm is compared with the optimal assignment method and with a greedy assignment method.

\end{enumerate}
Thus, our solution approach for the car parking problem is fair, distributed, and is easily deployed with the coordination of a central entity. In addition, it has appealing privacy properties.

The rest of the paper is organized as follows. A description of the system model and the problem formulation is presented in Section~\ref{sec:model}. In Section~\ref{sec:solving}, we provide the solution method to the car parking problem by using duality theory and subgradient method. Section~\ref{sec:ald_devlop} presents our proposed algorithm for distributed car parking assignment problem. In Section~\ref{sec:privacy}, we describe privacy properties of the algorithm. In Section~\ref{sec:results}, numerical results are provided. Lastly, Section~\ref{sec:conclusions} concludes the paper.

\subsection*{Notations}
Boldface lower case and upper case letters represent vectors and matrices, respectively, and calligraphy letters represent sets. The set of real $n$-vectors is denoted by $\R^n$ and the set of real $m\times n$ matrices is denoted $\R^{m\times n}$. We use {parentheses} to construct matrices from comma separated sub-matrices of agreed dimensions, e.g., $({\vec A}, {\vec B}, {\vec C})=[{\vec A}\tran \ {\vec B}\tran \ {\vec C}\tran]\tran$. We denote by $({\vec A}_i)_{i=1,2,\ldots,N}$ the matrix ${\vec A}=({\vec A}_1,{\vec A}_2,\ldots, {\vec A}_N)$.
The cardinality of a set ${\mathcal A}$ is denoted by $\mathrm{card} \ {\mathcal A}$.

\begin{table*}[t]
\centering
\subtable[]{
\begin{tabular}{l|ccccc}
\hline
Car index & \multicolumn{5}{c}{Free Parking Slot, $j$} \\
\cline{2-6}
$i$ & $j=1$ & $2$ & $3$  & $4$  & $5$\\
\cline{1-6}
$1$ & $0$ &$1$  &$0$ &$0$ &$0$\\
$2$ & $0$ &$0$  &$1$ &$0$ &$0$\\
$3$ & $0$ &$1$  &$0$ &$0$ &$0$\\
\bottomrule
\end{tabular}
\label{tab:incorrect}}
\qquad \qquad
\subtable[]{
\begin{tabular}{l|ccccc}
\hline
Car index & \multicolumn{5}{c}{Free Parking Slot, $j$} \\
\cline{2-6}
$i$ & $j=1$ & $2$ & $3$  & $4$  & $5$\\
\cline{1-6}
$1$ & $0$ &$1$  &$0$ &$0$ &$0$\\
$2$ & $0$ &$0$  &$1$ &$0$ &$0$\\
$3$ & $1$ &$0$  &$0$ &$0$ &$0$\\
\bottomrule
\end{tabular}
\label{tab:correct}}
\caption{\small{Assignment: $\mathcal{N}_t=\{1,2,3\}$, $\mathcal{M}_t=\{1,2,3,4,5\}$: (a) An infeasible assignment; (b) A feasible assignment.}}
\label{tab:Table-check}
\vspace{-5mm}
\end{table*}

\section{System model and Problem Formulation}  \label{sec:model}

A system consisting of $M$ parking slots and a number of destinations is considered. We denote by $\mathcal{M}=\{1,\ldots,M\}$ the set of parking slots. Destinations can include \emph{any} geographical location, such as shops, bars, banks, cinemas, casinos, houses, parks, hotels, among others. The parking slots and the destinations can be geographically dispersed and need not necessarily be concentrated. Knowledge of geographical location of each parking slot is assumed to be available to anyone in the system. A \emph{trustworthy cental controller} (CC) is responsible for coordinating the parking slot assignment~mechanism, namely it is the central authority. The coordinations are carried out through secured channels.

The parking slot assignment mechanism is assumed to operate in slotted time, with the
slots normalized to integer values $t\in\{1, 2, 3, \ldots\}$. At the beginning of every time slot~$t$, the set $\mathcal{M}_t\subseteq\mathcal{M}$ of \emph{free} parking slots is known.~\footnote{Such information is retrieved by installing sensors at every parking slots.} In addition, at the beginning of every time slot~$t$, a set $\mathcal{N}_t=\{1,\ldots,N_t\}$ of cars is scheduled for parking slot assignment, where $N_t\leq |M_t|$ is the total number of cars.
We denote by $\mbox{des}(i)$ the destination of car~$i\in\mathcal{N}_t$ and by $d_{ij}$ the distance from $\mbox{des}(i)$ to free parking slot $j\in\mathcal{M}_t$. We assume that each car~$i$ can compute $\{d_{ij}\}_{j\in\mathcal{M}_t}$ simply by knowing the geographical location of $\mbox{des}(i)$. Such computations can easily be performed by using the state-of-the-art global positioning system (GPS).  

To formally express the problem, let us first introduce \emph{binary decision variables} $(x_{ij})_{i\in\mathcal{N}_t, \ j\in\mathcal{M}_t}$, which indicate the $i$ to $j$ assignments as follows:
\be\label{eq:decesion_variables}
x_{ij}= \left\{ \begin{array}{ll}
  1 & \ \ \textrm{car $i$ is assigned to parking slot $j$}\\
  0  & \ \ \mbox{otherwise} \ .
   \end{array} \right.
\ee
A \emph{feasible assignment} should be such that one car is assigned to only one free parking slot and no more than one car is assigned to a free parking slot. For example, Table~\ref{tab:incorrect} shows an infeasible assignment and Table~\ref{tab:correct} shows a feasible assignment. Now we can formally express the distance from car~$i$'s assigned parking slot to its destination $\mbox{des}(i)$ as $\sum_{j\in\mathcal{M}_t}d_{ij}x_{ij}$. We refer to $\sum_{j\in\mathcal{M}_t}d_{ij}x_{ij}$ as \emph{parking distance} of car~$i$.


In order to ensure {min-max fairness} among the cars, we \emph{minimize the maximum parking distance}. Min-max fairness is appealing in many application domains in the sense that it ensures equalization of the costs incurred by the users, see~\cite{Radunovic-Boudec-2007}. Specifically, the problem can be formally expressed~as
\begin{IEEEeqnarray}{lcl}\label{eq:primal}
\mbox{minimize} & \ \ & \displaystyle\mbox{max}_{i\in\mathcal{N}_t}\ \textstyle\sum_{j \in \mathcal{M}_t} d_{ij}\ x_{ij}\IEEEyessubnumber\label{eq:problem_formulation}\\
\mbox{subject to} & \ \  & \textstyle{\sum_{i\in\mathcal{N}_t}} x_{ij} \leq 1, \ j\in\mathcal{M}_t \IEEEyessubnumber\label{eq:problem_formulation1}\\
& \ \ & \textstyle{\sum_{j\in\mathcal{M}_t}} x_{ij} = 1, \   i\in\mathcal{N}_t \IEEEyessubnumber\label{eq:problem_formulation2} \\
& \ \ & x_{ij}\in\{0,1\}, \quad i\in\mathcal{N}_t, \ j\in\mathcal{M}_t \IEEEyessubnumber\label{eq:problem_formulation3} \ ,
\end{IEEEeqnarray}
where the variable is $(x_{ij})_{i\in\mathcal{N}_t, \ j\in\mathcal{M}_t}$. Constraint~(\ref{eq:problem_formulation1}) ensures that no more than one car is assigned to a free parking slot. Constraint~(\ref{eq:problem_formulation2}) imposes that each car is assigned to only one free slot. Finally, constraint~(\ref{eq:problem_formulation3}) ensures that the values of $x_{ij}$ are either $0$ or $1$.

Note that the problem is nonconvex and even combinatorial. Hence we have to rely on global optimal methods~\cite{12} such as exhaustive search and branch and bound methods to \emph{solve} it. The main disadvantage of global methods is the prohibitive computational complexity, even in the case of small problems. Such methods are not scalable, and therefore can be impractical. In the sequel, we provide a method based on duality. Even though the optimality cannot be guaranteed, the proposed method is efficient, fast, and allows distributed implementation with a little coordination from the CC.


\section{Solution approach via dual problem }\label{sec:solving}


In this section, we first equivalently formulate problem~(\ref{eq:primal}) in its its epigraph form~\cite{Boyd-Vandenberghe-04}. Then we apply duality theory to obtain the related dual problem, and show that the problem is split into subproblems and a master problem which can be solved efficiently.

The equivalent problem is given by~\footnote{Without loss of generality, we drop the subindex~$t$ for notational simplicity.}
%
%
\begin{IEEEeqnarray}{lcl}\label{eq:epigraph}
\mbox{minimize} & \ \ & s \IEEEyessubnumber\label{eq:epi1}\\
\mbox{subject to} & \ \ &  \textstyle\sum_{j \in \mathcal{M}} d_{ij}\ x_{ij} \leq s, \  i\in\mathcal{N}\ \IEEEyessubnumber\label{eq:epi2}\\
& \ \  & \textstyle{\sum_{i\in\mathcal{N}}} x_{ij} \leq 1, \  j\in\mathcal{M} \IEEEyessubnumber\label{eq:epi3}\\
& \ \ & \textstyle{\sum_{j\in\mathcal{M}}} x_{ij} = 1, \   i\in\mathcal{N} \IEEEyessubnumber\label{eq:epi4} \\
& \ \ & x_{ij}\in\{0,1\}, \  i\in\mathcal{N}, \ j\in\mathcal{M} \IEEEyessubnumber\label{eq:epi5} \ ,
\end{IEEEeqnarray}
where the variables are $s$ and ${\vec x} = (x_{ij})_{i\in\mathcal{N}, \ j\in\mathcal{M}}$. Note that like problem~(\ref{eq:primal}), \eqref{eq:epigraph} is still nonconvex. Now we seek to decouple the problem among the cars, in order to maintain scalability properties of the car parking mechanism. In this context, we can clearly see that constraints~\eqref{eq:epi4},~\eqref{eq:epi5} are already decoupled, yet constraints~\eqref{eq:epi2},~\eqref{eq:epi3} are coupled among the cars, which is an obstacle to distributed solution methods.

Let us now form the partial Lagrangian by dualizing the coupling constraints~\eqref{eq:epi2} and \eqref{eq:epi3}. To do this, we introduce multipliers $\boldsymbol \lambda = (\lambda_i)_{i\in\mathcal{N}}$ for the inequality constraints~\eqref{eq:epi2} and multipliers $\boldsymbol \mu = (\mu_j)_{j\in\mathcal{M}}$ for the inequality constraints~\eqref{eq:epi3}.
The Lagrangian associated with problem~(\ref{eq:epigraph}) is
\begin{equation}
\begin{split} \label{eq:lagrangian}
\hspace{-3mm}L(s,{\vec x},{\boldsymbol\lambda}, {\boldsymbol\mu}) & {=} s{+} \ssum{i \in \mathcal{N}} \lambda_i \Big( \ssum{j \in \mathcal{M}} d_{ij}x_{ij} {-} s\Big) \\
 & \hspace{30mm} + \ssum{j \in \mathcal{M}} \mu_j \Big(\ssum{i \in \mathcal{N}} x_{ij}{-}1\Big) \\
& = s\Big(1{-} \ssum{i \in \mathcal{N}} \lambda_i\Big) {+} \ssum{i \in \mathcal{N}} \ssum{j \in \mathcal{M}} (\lambda_i d_{ij} {+} \mu_j) x_{ij} \\
& \hspace{47mm} - \ssum{j \in \mathcal{M}} \mu_j \ .
\end{split}
\end{equation}
\\
The dual function $g(\boldsymbol\lambda,\boldsymbol\mu)$ is given by
\begin{subequations}\label{eq:dual_function}
\begin{align}\label{eq:dual_function_}
   \hspace{-1mm}& g\big( \boldsymbol{\lambda},\boldsymbol{\mu}\big) = \inf_{\substack{s \in \R, \\ \sum_{j \in \mathcal{M}} x_{ij}=1, \ i\in\mathcal{N}, \\  x_{ij} \in \left\{0,1 \right\}, \ i\in\mathcal{N}, \ j\in\mathcal{M} }}L\big(s,{\vec x},\boldsymbol{\lambda},\boldsymbol{\mu}\big)\\ \label{eq:dual_function_2} \displaybreak[1]
  & \hspace{-2mm}=\left\{\hspace{-3mm} \begin{array}{ll}
  \displaystyle\inf_{\substack{\sum_{j \in \mathcal{M}} x_{ij}=1, i\in\mathcal{N}, \\  x_{ij} \in \left\{0,1 \right\}, i\in\mathcal{N}, \ j\in\mathcal{M} }}\hspace{-1mm} \ssum{i \in \mathcal{N}} \ssum{j \in \mathcal{M}} (\lambda_i d_{ij} {+} \mu_j) x_{ij} {-}\hspace{-3mm} \ssum{j \in \mathcal{M}} \mu_j \\
  \hspace{52mm} \mathrm{if} \ \displaystyle\textrm{$\ssum{i \in \mathcal{N}} \lambda_{i}=1$} \\
  -\infty \qquad \mathrm{otherwise} 
   \end{array} \right.\\ \label{eq:dual_function_3} \displaybreak[1]
  & \hspace{-2mm}= \left\{\hspace{-3mm} \begin{array}{ll}
   \displaystyle\ssum{i \in \mathcal{N}}\hspace{-1mm} \Bigg( \displaystyle\inf_{\substack{\sum_{j \in \mathcal{M}} x_{ij}=1, \\ x_{ij} \in \left\{0,1 \right\}, \ j\in\mathcal{M} }} \hspace{-2mm}\ssum{j \in \mathcal{M}} (\lambda_i d_{ij} {+} \mu_j) x_{ij}\hspace{-1mm}\Bigg) {-} \ssum{j \in \mathcal{M}} \mu_j\\  \hspace{52mm} \mathrm{if} \ \displaystyle\textrm{$\ssum{i \in \mathcal{N}} \lambda_{i}=1$} \\
  -\infty  \qquad \mathrm{otherwise} 
   \end{array} \right.\\ \label{eq:dual_function_4} \displaybreak[1]
  &\hspace{-2mm}= \left\{ \begin{array}{ll}
  \displaystyle  \ssum{i \in \mathcal{N}} g_i(\boldsymbol{\lambda},\boldsymbol{\mu})  -\ssum{j \in \mathcal{M}}\mu_j& \qquad \displaystyle\textrm{$\mathop{\textstyle{\sum}}_{i \in \mathcal{N}} \lambda_{i}=1$}\\
  -\infty  & \qquad \mathrm{otherwise} \ ,
   \end{array} \right.
\end{align}
\end{subequations}
where the equality~(\ref{eq:dual_function_2}) follows from that the linear function $s(1-\sum_{i\in\mathcal{N}}\lambda_i)$ is bounded below only when it is identically zero, the equality~(\ref{eq:dual_function_3}) follows from that constraints~\eqref{eq:epi4},~\eqref{eq:epi5} are separable, and $g_i(\boldsymbol{\lambda},\boldsymbol{\mu})$ is the optimal value of the problem
\begin{equation} \label{eq:sub_problem}
\begin{array}{ll}
\mbox{minimize} &  \ssum{j \in \mathcal{M}} (\lambda_i d_{ij}+\mu_j) x_{ij}\\
\mbox{subject to} & \ssum{j \in \mathcal{M}} x_{ij}=1 \\
& x_{ij} \in \left\{0,1 \right\}, \ j\in\mathcal{M}  \ ,
\end{array}
\end{equation}
with the variable $(x_{ij})_{j\in\mathcal{M}}$. Note that problem~\eqref{eq:sub_problem} is a combinatorial problem. Nevertheless, it has a closed-form solution
\be\label{eq:sub_problem_soln}
x^\star_{ij}= \left\{ \begin{array}{ll}
  1 & \ \ \textrm{$j=\displaystyle\mathop{\arg\min}_{l\in\mathcal{M}}(\lambda_i d_{il}+\mu_l)$}\\
  0  & \ \ \mathrm{otherwise} \ .
   \end{array} \right.
\ee
The dual master problem is
\begin{IEEEeqnarray}{lcl}\label{eq:dual_problem}
\mbox{maximize} & \ \ & g(\boldsymbol{\lambda},\boldsymbol{\mu})=\textstyle\sum_{i \in \mathcal{N}} g_i(\boldsymbol{\lambda},\boldsymbol{\mu})\IEEEyessubnumber\label{eq:dual_problem1}\\
\mbox{subject to} & \ \  & \textstyle\sum_{i \in \mathcal{N}} \lambda_{i}=1 \IEEEyessubnumber\label{eq:dual_problem2}\\
& \ \ & \textstyle\lambda_i\geq 0, \ i\in\mathcal{N} \IEEEyessubnumber\label{eq:dual_problem3} \\
& \ \ & \textstyle\mu_j\geq 0, \ j\in\mathcal{M} \IEEEyessubnumber\label{eq:dual_problem4}  \ ,
\end{IEEEeqnarray}
where the variables are ${\boldsymbol\lambda}$ and ${\boldsymbol\mu}$. In the sequel, we describe an approach to solve the dual problem~\eqref{eq:dual_problem}, based on the projected subgradient method~\cite{Boyd-EE364b-SubGradMethods-07}.

\subsection*{Solving the dual}
Note that $g(\boldsymbol\lambda, \boldsymbol\mu)$ is a concave function, therefore,  we need to find the subgradient ${\vec s}\in\R^{N+M}$ of $-g$ at a feasible $(\boldsymbol\lambda, \boldsymbol\mu)$. For clarity we separate ${\vec s}$  into two vectors as follows:
\begin{equation}
{\vec s}=({\vec u},{\vec v}),
\end{equation}
where ${\vec u}=(u_i)_{i\in\mathcal{N}}$ is the part of ${\vec s}$ that corresponds to $\boldsymbol\lambda$ and ${\vec v}=(v_j)_{j\in\mathcal{M}}$ the part that corresponds to $\boldsymbol\mu$. The negative of dual function $-g(\boldsymbol\lambda,\boldsymbol\mu)$ is given by
\begin{equation*}
\begin{split}
 -g(\boldsymbol\lambda, \boldsymbol\mu) & = \ssum{j \in \mathcal{M}} \mu_j - \ssum{j \in \mathcal{M}} \mu_j \ssum{i \in \mathcal{N}} x^{\star}_{ij} - \ssum{i \in \mathcal{N}} \lambda_i \ssum{j \in \mathcal{M}} d_{ij}x^{\star}_{ij} \ ,
\end{split}
\end{equation*}
and particular choices for $u_i$, $i \in \mathcal{N}$ and $v_j$, $j \in \mathcal{M}$ are given by
\begin{equation} \label{eq:uv}
u_i = - \ssum{j \in \mathcal{M}} d_{ij} x^{\star}_{ij} \qquad \mbox{and} \qquad v_j= 1- \ssum{i \in \mathcal{N}} x^{\star}_{ij} \ .
\end{equation}
Thus the projected subgradient method is given by
\begin{equation}  \label{eq:it}
({\boldsymbol\lambda}^{(k+1)}, {\boldsymbol\mu}^{(k+1)}) = P(({\boldsymbol\lambda}^{(k)}, {\boldsymbol\mu}^{(k)}) - \alpha_k ({\vec u}^{(k)},{\vec v}^{(k)})) \ ,
\end{equation}
where $k$ is the current iteration index of the subgradient method, $P({\vec z})$ is the Euclidean projection of ${\vec z}\in\R^{N+M}$ onto the \emph{feasible set} of the dual problem~\eqref{eq:dual_problem}, and $\alpha_k > 0$ is the $k$th step size, chosen to guarantee the asymptotic convergence of the subgradient method, e.g., $\alpha_k = \alpha / k$, where $\alpha$ is a positive scalar. Since the feasible set of dual problem is separable in $\lambda$ and $\mu$, the projection $P( \cdot)$  can be performed independently. Therefore, the iteration~(\ref{eq:it}) is equivalently performed as follows:
\begin{equation}   \label{eq:it1}
{\boldsymbol\lambda}^{(k+1)} = P_s({\boldsymbol\lambda}^{(k)} - \alpha_k {\vec u}^{(k)})
\end{equation}
\begin{equation}   \label{eq:it2}
{\boldsymbol\mu}^{(k+1)} = [{\boldsymbol\mu}^{(k)} - \alpha_k {\vec v}^{(k)}]^+,
\end{equation}
where $P_s( \cdot)$ is the Euclidean projection onto the probability simplex~\cite{Bertsekas-99},
\begin{equation}  \label{eq:simplex}
\Pi = \left\{ \lambda \left| \textstyle\sum^N_{i=1} \lambda_i =1, \lambda_i \geq 0 \right. \right\}
\end{equation}
and $ [ \  \cdot \  ]^+$ is the Euclidean projection onto the nonnegative orthant. Note that the Euclidean projection onto the probability simplex can be posed as a convex optimization problem that can be solved efficiently, see Appendix~\ref{app:ProjSimplex}.

\section{Algorithm Development}\label{sec:ald_devlop}
In this section, we first present our distributed algorithm to address problem~\eqref{eq:primal} via the dual problem~\eqref{eq:dual_problem}. The resulting algorithm is indeed the distributed car parking mechanism that can be coordinated by the CC or the central controller. Next, we discuss the convergence properties of the algorithm.

\subsection{Distributed algorithm implementation }\label{subsec:DPP}
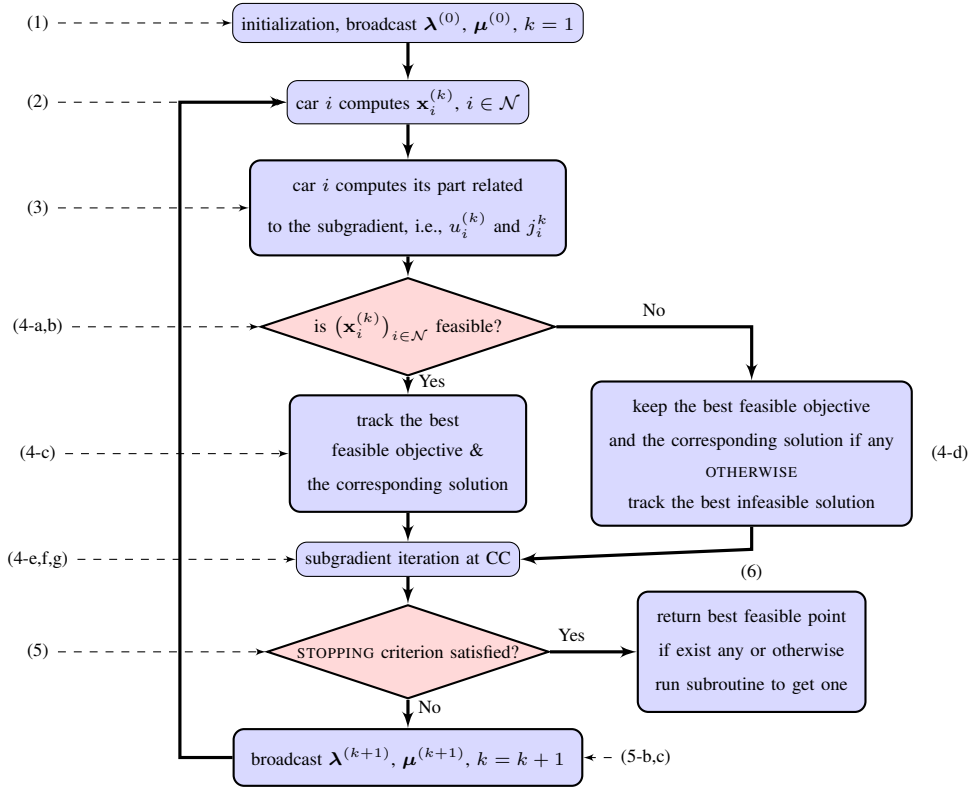
\begin{figure*}[!t]
\centering
\tikzstyle{decision} = [diamond, aspect=3, draw, thick, fill=red!15,minimum height=3mm,  text badly centered, inner sep=1pt,font=\scriptsize,node distance=4em]
\tikzstyle{MyMatrix1} = [matrix, draw, thick, fill=blue!15, text centered, rounded corners,font=\scriptsize,minimum height=0.5em,node distance=4em]
\tikzstyle{block} = [rectangle, draw, fill=blue!15,
     text centered, rounded corners, minimum height=1em,font=\scriptsize,node distance=3em]
\tikzstyle{line} = [draw, -latex']
\tikzstyle{cloud} = [draw=white, node distance=14em,font=\scriptsize]

\begin{tikzpicture}[scale=0.10]
    \node [block] (init) {initialization, broadcast $\boldsymbol{\lambda}^{(0)}$, $\boldsymbol{\mu}^{(0)}$, $k=1$};
    \node [cloud, left of=init] (step1) {(1)};

    \node [block, below of=init] (subprob) {car~$i$ computes ${\vec x}^{(k)}_i$, $i\in\mathcal{N}$};
    \node [cloud, left of=subprob] (step2) {(2)};

    \node [MyMatrix1,below of=subprob] (subgrad)
    {
        \node  {car~$i$ computes its part related}; \\
        \node  {to the subgradient, i.e., $u^{(k)}_i$ and $j^k_i$}; \\
    };
    \node [cloud, left of=subgrad] (step3) {(3)};

    \node [decision, below of=subgrad,node distance=4.5em] (feasibility) {is $\big({\vec x}^{(k)}_i\big)_{i\in\mathcal{N}}$ feasible?};
    \node [cloud, left of=feasibility] (step4a) {(4-a,b)};

    \node [MyMatrix1,below of=feasibility,node distance=4.8em] (bestFeObj)
    {
        \node  {track the best}; \\
        \node  {feasible objective \&}; \\
        \node {the corresponding solution };\\
    };
    \node [cloud, left of=bestFeObj] (step4c) {(4-c)};

    \node [MyMatrix1,right of=bestFeObj,node distance=13em] (bestInfObj)
    {
        \node  {keep the best feasible objective}; \\
        \node  {and the corresponding solution if any};\\
        \node  {\textsc{otherwise}}; \\
        \node  {track the best infeasible solution}; \\
    };
    \node [cloud, right of=bestInfObj,node distance=7.5em] (step4d) {(4-d)};

    \node [block, below of=bestFeObj,node distance=4.0em] (subgradIter) {subgradient iteration at CC};
    \node [cloud, left of=subgradIter] (step4efg) {(4-e,f,g)};

    \node [decision, below of=subgradIter,node distance=3.5em] (stop) {\textsc{stopping} criterion satisfied?};
    \node [cloud, left of=stop] (step5) {(5)};

    \node [MyMatrix1,right of=stop,node distance=13em] (suroutine)
    {
        \node  {return best feasible point}; \\
        \node  {if exist any or otherwise }; \\
        \node  { run subroutine to get one}; \\
    };
    \node [cloud, above of=suroutine,node distance=3.0em] (step6) {(6)};

    \node [MyMatrix1,below of=stop] (comPrices)
    {
        \node  {broadcast $\boldsymbol{\lambda}^{(k+1)}$, $\boldsymbol{\mu}^{(k+1)}$, $k=k+1$}; \\
    };
    \node [cloud, right of=comPrices,node distance=9.0em] (step5cd) {(5-b,c)};


    \path [line, very thick] (init) -- (subprob);
    \path [line, very thick] (subprob) -- (subgrad);
    \path [line, very thick] (subgrad) -- (feasibility);
    \path [line, very thick] (feasibility) -- node [near start, right] {\scriptsize{Yes}} (bestFeObj);
    \path [line, very thick] (bestFeObj) -- (subgradIter);
    \path [line, very thick] (subgradIter) -- (stop);
    \path [line, very thick] (stop) -- node [near start, right] {\scriptsize{No}} (comPrices);
    \path [line, very thick] (comPrices.west) -- ++(-7,0)|- (subprob);
    \path [line, very thick] (feasibility) -| node [near start, above] {\scriptsize{No}} (bestInfObj);
    \path [line, very thick] (bestInfObj.south) -- ++(0,-3)--  (subgradIter.east);
    \path [line, very thick] (stop) -- node [near start, above] {\scriptsize{Yes}} (suroutine);



    \path [line,dashed] (step1) -- (init);
    \path [line,dashed] (step2) -- (subprob);
    \path [line,dashed] (step3) -- (subgrad);
    \path [line,dashed] (step4a) -- (feasibility);
    \path [line,dashed] (step4c) -- (bestFeObj);
    \path [line,dashed] (step4efg) -- (subgradIter);
    \path [line,dashed] (step5) -- (stop);
    \path [line,dashed] (step5cd) -- (comPrices);
\end{tikzpicture}
\captionof{figure}{\label{fig:DAA}\small{Block diagram of the DCP algorithm.}}
\vspace{-5mm}
\end{figure*}

Roughly speaking, the algorithm  capitalizes on the ability of the CC to construct the subgradient $({\vec u},{\vec v})$ in a distributed fashion via the coordination of scheduled cars. Note that, the involvement of a CC (e.g., an authority who handles the parking slots) is essential for realizing the overall algorithm in practice. This involvement is mainly for coordinating certain parameter among the scheduled cars, and for constricting a feasible assignment in case the assignment from dual problem is infeasible. The algorithm is formally documented below, see also Fig.~\ref{fig:DAA} for a concise block diagram.

\noindent\rule{0.49\textwidth}{0.3mm}
\\
\emph{Algorithm}: \ \textsc{\small{Distributed car-parking (\textrm{DCP})}}
\begin{enumerate}
\item Given the distances $(d_{ij})_{j\in\mathcal{M}}$ for each car~$i\in\mathcal{N}$. The central controller (CC) sets $k=1$, sets current objective value $p^{\textrm{cur}}(0)=\infty$, sets number of conflicting users $N^{\textrm{conflict}}=N$, and broadcasts the initial (feasible) $\lambda^{(k)}_i$ and $\boldsymbol\mu^{(k)}=(\mu^{(k)}_j)_{j\in\mathcal{M}}$ to each car~$i\in\mathcal{N}$.
\item Every car $i$ sets $\lambda_i = \lambda_i^{(k)}$ and $\boldsymbol\mu=\boldsymbol\mu^{(k)}$ and locally computes ${\vec x}^{(k)}_i=(x^\star_{ij})_{j\in\mathcal{M}}$ from~\eqref{eq:sub_problem_soln}. Let $j^{k}_i$ denote the index of the nonzero component of ${\vec x}^{(k)}_i$.
\item Local subgradients: Each car $i$ 
	\begin{itemize}
			\item[a.] sets \emph{scalar} $u^{(k)}_i=- \ssum{j \in \mathcal{M}} d_{ij} x^{(k)}_{ij}= -d_{ij^{k}_i}$, [compare with~(\ref{eq:uv})].
            \item[b.] transmits $(u^{(k)}_i,j^{k}_i)$ to CC.
		\end{itemize}
\item Current assignment and Subgradient iteration at CC
		\begin{itemize}
            \item[a.] find the set $\mathcal{J}^{(k)}_j$ of users assigned to slot~$j$, i.e., $\mathcal{J}^{(k)}_j=\{i \ | \ j^k_i=j\}$. Set $N^{\textrm{conflict}}_k=\sum_{j|\mathrm{card}(\mathcal{J}^{(k)}_j)\geq 2}\mathrm{card}(\mathcal{J}^{(k)}_j)$.
			\item[b.] if no conflicting assignments (i.e., $N^{\textrm{conflict}}_k=0$), set $N^{\textrm{conflict}}=N^{\textrm{conflict}}_k$ and go to step~4-c. Otherwise, go to step~4-d.
    \item[c.] if $p^{\textrm{cur}}(k-1)>\max_{i\in\mathcal{N}}d_{ij^{k}_i}$, set $p^{\textrm{cur}}(k)=\max_{i\in\mathcal{N}}d_{ij^{k}_i}$ and set ${\vec X}^{\textrm{cur}}(k)=\big({\vec e}\tran_{j^{k}_i}\big)_{i\in\mathcal{N}}\in\R^{N\times M}$. Go to step~4-e.
    \item[d.] if $N^{\textrm{conflict}}_k< N^{\textrm{conflict}}$, set $N^{\textrm{conflict}}=N^{\textrm{conflict}}_k$, $p^{\textrm{cur}}(k)=\infty$, ${\vec X}^{\textrm{cur}}(k)=\big({\vec e}\tran_{j^{k}_i}\big)_{i\in\mathcal{N}}\in\R^{N\times M}$, and $\mathcal{J}^{\textrm{cur}}_j=\mathcal{J}^{(k)}_j$, ${j\in\mathcal{M}}$. Go to step~4-e.

			\item[e.] form ${\vec u}^{(k)}=(u^{(k)}_i)_{i\in\mathcal{N}}$ and perform~(\ref{eq:it1}) to find $\boldsymbol\lambda^{(k+1)}$.
			\item[f.] set $v_j= 1-\mathrm{card}(\mathcal{J}_j)$, ${\vec v}^{(k)}=(v_j)_{j\in \mathcal{M}}$, [compare with~(\ref{eq:uv})].
            \item[g.] perform~(\ref{eq:it2}) to find $\boldsymbol\mu^{(k+1)}$.
		\end{itemize}
\item Stopping criterion: If the stopping criterion is satisfied,
    \begin{itemize}
			\item[a.] go to step~6.
		\end{itemize}
Otherwise,
        \begin{itemize}
        \item[b.] CC broadcasts the new $\lambda^{(k+1)}_i$ and $\boldsymbol\mu^{(k+1)}$ to each car~$i\in\mathcal{N}$.
        \item[c.] increment $k$, i.e., set $k=k+1$.
		\item[d.] go to step 2.
		\end{itemize}
\item Output: If $N^{\textrm{conflict}}=0$ (i.e., a feasible assignment is achieved), CC returns ${\vec X}^{\textrm{final}}={\vec X}^{\textrm{cur}}(k)$ and terminates the algorithm. Otherwise, CC sets ${\vec X}^{\textrm{infeasible}}={\vec X}^{\textrm{cur}}(k)$, performs a simple \emph{subroutine} to construct a feasible assignment ${\vec X}^{\textrm{final}}$ from ${\vec X}^{\textrm{infeasible}}$, returns ${\vec X}^{\textrm{final}}$, and terminates the algorithm.
\end{enumerate}
\vspace{-3mm}
\rule{0.49\textwidth}{0.3mm}\vspace{-5mm}
\subsection{Algorithm description}
In step~1, the algorithm starts by choosing initial feasible values for $\lambda^{(k)}_i, \ i\in\mathcal{N}$ and $\mu^{(k)}_j, \ j\in\mathcal{M}$.
Step~2 corresponds to the local computations of ${\vec x}^{(k)}_i$ at each car~$i$. These computations involve simple comparisons [see~\eqref{eq:sub_problem_soln}] and can be performed in a parallelized manner by the scheduled~cars.
Step~3 involves coordination of scheduled cars and CC. First, each car~$i$ constructs \textit{scalar} parameter $u^{(k)}_i$. Then it transmits $u^{(k)}_i$ together with the potential car slot index $j^{k}_i$ to CC. 


In step~4, CC keeps records of the \emph{best} assignment so far. The assignment is best, in the following sense. First suppose algorithm yields at least one feasible assignment, i.e., $N^{\textrm{conflict}}=0$. Then the best assignment is the one that corresponds to the smallest objective value among all feasible assignments, see step~4-c. On the other hand, suppose algorithm does not yield any feasible assignment, i.e., $N^{\textrm{conflict}}>0$. Then the best assignment is the one corresponds to the smallest $N^{\textrm{conflict}}_k$ among all infeasible assignments, see step~4-d. Note that $N^{\textrm{conflict}}$ is equal to the total conflicting users, and thus quantifies the degrees of infeasibility, see step~4-a.
Moreover, by using the information received from scheduled cars, CC constructs the global subgradient components ${\vec u}^{(k)}\in\R^N$ and ${\vec v}^{(k)}\in\R^M$, which in turn are used to perform the subgradient iterations~\eqref{eq:it1}-\eqref{eq:it2}, see steps~4-e,4-f,4-g.

The new parameters $\boldsymbol \lambda^{(k+1)}$ and $\boldsymbol \mu^{(k+1)}$ are broadcasted to every car and the algorithm is performed in an iterative manner until a stopping criterion is~satisfied, see step~5. Natural stopping criterion in practice includes running subgradient method for a fixed number of iterations.

Recall that, solution method for the primal problem~\eqref{eq:epigraph} by considering its dual problem~\eqref{eq:dual_problem} does not always guarantee the primal feasibility, because the original problem~\eqref{eq:epigraph} is \emph{nonconvex}~\cite{Bertsekas-99}. Therefore, if a feasible assignment is not achieved, a subroutine call is required to construct one after the stopping criterion is satisfied. Step~6 is essentially to address this infeasibility problem. In particular, once the stopping criterion is satisfied (step~5), CC checks whether the current assignment ${\vec X}^{\textrm{cur}}(k)$ obtained is feasible. If it is feasible, algorithm terminates by returning ${\vec X}^{\textrm{final}}={\vec X}^{\textrm{cur}}(k)$, where CC informs each car~$i$, its parking slot. Otherwise, CC performs a simple subroutine to construct a feasible assignment ${\vec X}^{\textrm{final}}$ by using the current infeasible assignment~${\vec X}^{\textrm{infeasible}}$, before the algorithm terminates. In the sequel, we outline a subroutine that can be implemented at CC for constructing a feasible assignment.

\subsection{Constructing a feasible assignment from ${\vec X}^{\textrm{infeasible}}$}\label{subsec:subroutine}
The key idea of the subroutine is summarized as follows: 1) select the set of cars that are assigned to the \emph{same} parking slot, 2) find the set of free parking slots, and 3) assign the conflicting cars found in the first stage to the free parking slots found in the second stage in an iterative manner. In the following, we describe this idea in the detail.

We start by introducing some useful notations for clarity. We denote by $\mathcal{M}^{\textrm{over-assigned}}$ the set of parking slots, where two or more than two cars are assigned, i.e., $\mathcal{M}^{\textrm{over-assigned}}= \{j \ | \ \mathrm{card}(\mathcal{J}^{\textrm{cur}}_j)\geq2 \}$. Moreover, we denote by $\mathcal{M}^{\textrm{free}}$ the set of free parking slots, i.e., $\mathcal{M}^{\textrm{free}}= \{ j\in\mathcal{M} \ | \ \mathcal{J}^{\textrm{cur}}_j=\emptyset  \}$.
For example, suppose Table~\ref{tab:incorrect} corresponds to the current assignment~${\vec X}^{\textrm{infeasible}}$, which is infeasible. Then we have $\mathcal{M}^{\textrm{over-assigned}}=\{2\}$ and $\mathcal{M}^{\textrm{free}}=\{1,4,5\}$.
%
To formally express the subroutine, it is further useful to introduce some minor notations, where we relabel the indices of cars and parking slots.
Let $\boldsymbol\sigma=(\sigma_l)_{l=1,\ldots,\mathrm{card}(\mathcal{M}^{\textrm{over-assigned}})}$ denote the parking slot indices~$j\in\mathcal{M}^{\textrm{over-assigned}}$ arranged in an increasing order. Moreover, we denote by $n_j$ the total cars assigned to $j$th parking slot \big[i.e., $\mathrm{card}(\mathcal{J}^{\textrm{cur}}_j)$\big], where $j\in\mathcal{M}^{\textrm{over-assigned}}$. Now, the subroutine can be formally expressed as follows:


%


\noindent\rule{0.49\textwidth}{0.3mm}
\\
\emph{Algorithm}: \ \textsc{\small{Construct a feasible assignment from}} ${\vec X}^{\textrm{infeasible}}$
\begin{enumerate}
\item Given the infeasible assignment ${\vec X}^{\textrm{infeasible}}$; $\mathcal{M}^{\textrm{over-assigned}}$; $\mathcal{M}^{\textrm{free}}$, $\boldsymbol\sigma$, and $n_j$ $\forall j\in\mathcal{M}^{\textrm{over-assigned}}$. Set ${\vec X}^{\textrm{final}}={\vec X}^{\textrm{infeasible}}$, $k=1$, and $l=1$.
\item CC sets $\boldsymbol\pi=(\pi_n)_{n=1,\ldots,n_{\sigma_l}}$ to be the car indices~$i\in\mathcal{J}^{\textrm{cur}}_{\sigma_l}$ arranged in an increasing order.
\item For $n=2:n_{\sigma_l}$
        \begin{itemize}
        \item[a.] CC sends $\mathcal{M}^{\textrm{free}}$ to car~$\pi_n$.
        \item[b.] car~$\pi_n$ chooses slot $j$, where $j=\arg\min_j{d_{\pi_nj}}$ and sends $j$ to CC.
        \item[c.] CC updates $\mathcal{M}^{\textrm{free}}$ as $\mathcal{M}^{\textrm{free}}=\mathcal{M}^{\textrm{free}}\setminus\{j\}$ and sets $[{\vec X}^{\textrm{final}}]_{\pi_nj}=1$.
		\end{itemize}

\item If $l=\mathrm{card}(\mathcal{M}^{\textrm{over-assigned}})$, return ${\vec X}^{\textrm{final}}$ and STOP. Otherwise, set $l=l+1$ and go to step~2.
\end{enumerate}
\vspace{-2mm}
\rule{0.49\textwidth}{0.3mm}\vspace{-0mm}
Step~1 is the initialization of the subroutine. Step~2 takes every over-assigned parking slots in the order defined by $\boldsymbol \sigma$ and the corresponding conflicting car indices are ordered as $\boldsymbol \pi$. In step~3, all of these cars, but $\pi_1$ are assigned to free parking slots in an iterative manner. In particular, the conflicting cars (except $\pi_1$) in the given over-assigned parking slot are given the free parking slot indices (see step~3-a) and every car chooses its parking slots in a greedy manner (see step~3-b). Note that the assignment of car~$\pi_1$ is not changed and it can remain in the slot already specified by ${\vec X}^{\textrm{final}}$, because the other cars are reassigned in step~3-a and step~3-b. Moreover, CC updates the assignment ${\vec X}^{\textrm{final}}$ accordingly, see step~3-c. Stopping criterion in step~4 checks whether all the cars in the over-assigned parking slots have been reassigned. If so, subroutine terminates by returning the feasible assignment ${\vec X}^{\textrm{final}}$. Otherwise, the subroutine continues by moving to the next over-assign parking slot. If the subroutine above is applied
to ${\vec X}^{\textrm{infeasible}}$ given in Table~\ref{tab:incorrect}, then a possible feasible assignment ${\vec X}^{\textrm{final}}$ is shown in Table~\ref{tab:correct}.

\subsection{Convergence}\label{sec:convergence}
In this section, we present the convergence properties of the proposed \textrm{DCP} algorithm for car parking. In particular, we show that for a sufficiently large number of subgradient iterations, the \textrm{DCP} algorithm converges to the dual optimal value of problem~\eqref{eq:dual_problem}. The convergence is established by the following proposition:

\textit{Proposition 1:} Let $g^{(k)}_{\textrm{best}} = \max \{  g(\boldsymbol\lambda^{(1)},\boldsymbol\mu^{(1)}),\ldots, g(\boldsymbol\lambda^{(k)},\boldsymbol\mu^{(k)})\}$ denote the dual objective value found after $k$ subgradient iterations and $(\boldsymbol\lambda^{\star}, \boldsymbol\mu^{\star})$ denote the optimal solution of dual problem~(\ref{eq:dual_problem}). Suppose $\|(\boldsymbol\lambda^{(1)},\boldsymbol\mu^{(1)})-(\boldsymbol\lambda^{\star}, \boldsymbol\mu^{\star})\|$ is bounded from above. Then, $ \forall \varepsilon >0, \  \exists n \geq 1 $ such that $ \forall k \  k\geq n \Rightarrow (d^\star - g^{(k)}_{\textrm{best}}) < \varepsilon$, where $d^\star$ is the optimal value of the dual problem~\eqref{eq:dual_problem}.
\begin{IEEEproof}
The proof is based on the approach of~\cite{Boyd-EE364b-SubGradMethods-07,Bertsekas-99}. We have
\begin{align}
&\hspace{-3mm}\big\|
(\boldsymbol\lambda^{(k+1)}, \boldsymbol\mu^{(k+1)})
-
(\boldsymbol\lambda^\star, \boldsymbol\mu^\star)\big\|^2_2 \displaybreak[1]\nonumber\\
& \hspace{-3mm}= \ \big\| P
\big((\boldsymbol\lambda^{(k)} {-} \alpha_k {\vec {\vec u}}^{(k)}){-}\boldsymbol\lambda^\star,(\boldsymbol\mu^{(k)} {-} \alpha_k {\vec v}^{(k)}){-}\boldsymbol\mu^\star\big)
\big\|^2_2 \label{eq:c1} \\ \displaybreak[1]
& \hspace{-3mm}\leq \ \big\|
((\boldsymbol\lambda^{(k)} - \alpha_k {\vec u}^{(k)})-\boldsymbol\lambda^\star, (\boldsymbol\mu^{(k)} - \alpha_k {\vec v}^{(k)})-\boldsymbol\mu^\star)
\big\|^2_2 \label{eq:c2} \\ \displaybreak[1]
& \hspace{-3mm}= \ \| (\boldsymbol\lambda^{(k)},\boldsymbol\mu^{(k)}) - (\boldsymbol\lambda^\star,\boldsymbol\mu^\star) \|^2_2 -2 \alpha_k {\vec u}^{(k)T}(\boldsymbol\lambda^{(k)} - \boldsymbol\lambda^{\star})\notag  \\ \displaybreak[1]
& \ \ \ -2 \alpha_k {\vec v}^{(k)T}(\boldsymbol\mu^{(k)} - \boldsymbol\mu^{\star}) +  \alpha^2_k \|{\vec u}^{(k)}\|^2_2 + \alpha^2_k \| {\vec v}^{(k)} \|^2_2  \label{eq:c3} \\ \displaybreak[1]
& \hspace{-3mm}\leq \ \| (\boldsymbol\lambda^{(k)},\boldsymbol\mu^{(k)}) - (\boldsymbol\lambda^\star,\boldsymbol\mu^\star) \|^2_2 -2 \alpha_k \big(g(\boldsymbol\lambda^\star,\boldsymbol\mu^{\star})  \nonumber \\ \displaybreak[1]
& \hspace{12mm}-g(\boldsymbol\lambda^{(k)},\boldsymbol\mu^{(k)})\big) + \alpha^2_k \|{\vec u}^{(k)} \|^2_2  + \alpha^2_k \|{\vec v}^{(k)} \|^2_2 \ , \label{eq:c5}
\end{align}
where~(\ref{eq:c1}) follows from~(\ref{eq:it}),~(\ref{eq:c2}) follows from that the Euclidean projection $P({\vec z})$ of any ${\vec z}\in\R^{N+M}$ onto the \emph{feasible set} of the dual problem~\eqref{eq:dual_problem} always decreases the distance of $P({\vec z})$ to every point in the \emph{feasible set} and in particular to the optimal point $(\boldsymbol\lambda^{\star}, \boldsymbol\mu^{\star})$, and~(\ref{eq:c5}) follows from the definition of subgradient. Recursively applying~(\ref{eq:c5}) and rearranging the terms we obtain
\begin{align}
& \hspace{-0mm}2 \textstyle{\ssum{l=1:k}} \alpha_l (d^{\star} - g(\boldsymbol\lambda^{(l)}, \boldsymbol\mu^{(l)}){\leq} {-} \big\|(\boldsymbol\lambda^{(k+1)},\boldsymbol\mu^{(k+1)})
{-}
(\boldsymbol\lambda^\star, \boldsymbol\mu^\star)\big\|^2_2 \nonumber\\\displaybreak[1]
& \hspace{-1mm}{+} \|(\boldsymbol\lambda^{(1)},\boldsymbol\mu^{(1)}){-}(\boldsymbol\lambda^{\star}, \boldsymbol\mu^{\star})\|^2_2   {+}\ssum{l=1:k} \alpha^2_l \| {\vec u}^{(l)}  \|^2_2 {+} \ssum{l=1:k} \alpha^2_l \| {\vec v}^{(l)}  \|^2_2 \label{eq:c6}\\\displaybreak[1]
& \hspace{-1mm}\leq \|(\boldsymbol\lambda^{(1)},\boldsymbol\mu^{(1)}){-}(\boldsymbol\lambda^{\star}, \boldsymbol\mu^{\star})\|^2_2 {+} \ssum{l=1:k} \alpha^2_l\big( \| {\vec u}^{(l)}  \|^2_2 {+} \| {\vec v}^{(l)}  \|^2_2\big)  \label{eq:c6_2} \\\displaybreak[1]
& \hspace{-1mm}\leq R^2 {+} G^2_1 \ssum{l=1:k} \alpha^2_l {+} G^2_2 \ssum{l=1:k} \alpha^2_l=R^2 {+} (G^2_1{+}G^2_2) \ssum{l=1:k} \alpha^2_l \ , \label{eq:c7}
\end{align}
where~(\ref{eq:c6_2}) follows from that $\big\|(\boldsymbol\lambda^{(k+1)},\boldsymbol\mu^{(k+1)})
-
(\boldsymbol\lambda^\star, \boldsymbol\mu^\star)\big\|^2_2  \geq 0$,
(\ref{eq:c7}) follows from that $\|(\boldsymbol\lambda^{(1)},\boldsymbol\mu^{(1)})-(\boldsymbol\lambda^{\star}, \boldsymbol\mu^{\star})\|$ is bounded from above, i.e., $\exists \ R<\infty$ such that $\|(\boldsymbol\lambda^{(1)},\boldsymbol\mu^{(1)})-(\boldsymbol\lambda^{\star}, \boldsymbol\mu^{\star})\|<R$ and the norm of the subgradient $({\vec u},{\vec v})$ is bounded from above as
\begin{align} \label{eq:c8}
\|{\vec u} \|_2 & \leq G_1 = \sqrt{\textstyle{\sum_{i \in \mathcal{N}}}(\textstyle{\max_{j \in \mathcal{M}}}{d_{ij}})^2} \\ \label{eq:c9}
\|{\vec v} \|_2 &  \leq G_2 = \sqrt{(N-1)^2+(M-1)} \ .
\end{align}
The bound  \eqref{eq:c8} is obtained by noting that there exist at most one nonzero element in $(x_{ij})_{j\in\mathcal{M}}$, see~\eqref{eq:sub_problem_soln} and \eqref{eq:uv}. Moreover, \eqref{eq:c9} follows when all cars are assigned to one parking slot, see \eqref{eq:uv}. By using the trivial relation
\begin{equation} \label{eq:c10}
d^{\star}-g^{(k)}_{\textrm{best}}  \leq d^\star - g(\boldsymbol\lambda^{(l)},\boldsymbol\mu^{(l)}), \ \ l = 1,...,k \ ,
\end{equation}
and \eqref{eq:c7},  we obtain an upper bound on $d^{\star}-g^{(k)}_{\textrm{best}}$ as
\begin{align} \label{eq:c11}
d^{\star}-g^{(k)}_{\textrm{best}} \ \leq \ & \big(R^2 +(G^2_1+G^2_2) \textstyle{\sum^k_{l=1}} \alpha^2_l\big) / (2 \textstyle{\sum^k_{l=1}} \alpha_l)
\end{align}
Noting that step size $\alpha_l=\alpha/l, \ 0<\alpha<\infty$ is \emph{square summable}, i.e., $\sum_{l=1}^\infty\alpha^2_l=\alpha^2\pi/6$. Moreover, $\sum_{l=1}^{k}\alpha_l$ is strictly monotonically increasing in $k$ (it grows without bound as $k\rightarrow\infty$). Therefore, for any $\epsilon >0$, we can always find an integer $n\geq 1$ such that $\sum_{l=1}^k\alpha_l> \epsilon \ (R^2/2 +a^2(G^2_1+G^2_2) \pi/12)$ if $k\geq n$, which concludes the proof.
\end{IEEEproof}
	
The bound derived in~\eqref{eq:c11}, together with \eqref{eq:c8}-\eqref{eq:c9} allows us to predict some key behaviors of the convergence of the proposed algorithm. For example, the larger the $d_{ij}$ values, the larger the $G_1$, and therefore, the larger the number of iterations to achieve a given accuracy. Nevertheless, the influence of $G_1$ can be made negligible by arbitrarily scaling down the objective function of problem~\eqref{eq:primal}. From~\eqref{eq:c9}, we note that the number of scheduled cars (i.e., $N$) and the number of free parking slots (i.e., $M$) directly influence the convergence. In the next section, we highlight some appealing privacy preserving properties of the \textrm{DCP} algorithm.	

\section{Privacy Properties of the algorithm} \label{sec:privacy}
We see that the proposed car parking mechanism DCP is preserving {privacy} in the sense that any car $n\neq i$ will not be able to find out the destination $\mbox{des}(i)$ of $i$th car while using the \textrm{DCP} algorithm. We refer to an attempt of an arbitrary car~$n$ to discover the destination of any other car $i$, as a passive attack~\cite[\S~5.1-5.3]{Goldreich-book-2004}, where car~$n$ keeps records of possibly all the information that it exchanges with CC and by using those it tries to discover private~data $\mbox{des}(i)$. In what follows, we first present sufficient information that an arbitrary car~$n$ can use to discover $\mbox{des}(i)$. Then we show how \textrm{DCP} algorithm cancels such a sufficient information and ensures privacy.

\subsection{Sufficient information to discover the destination} \label{subsec:tent_attack}
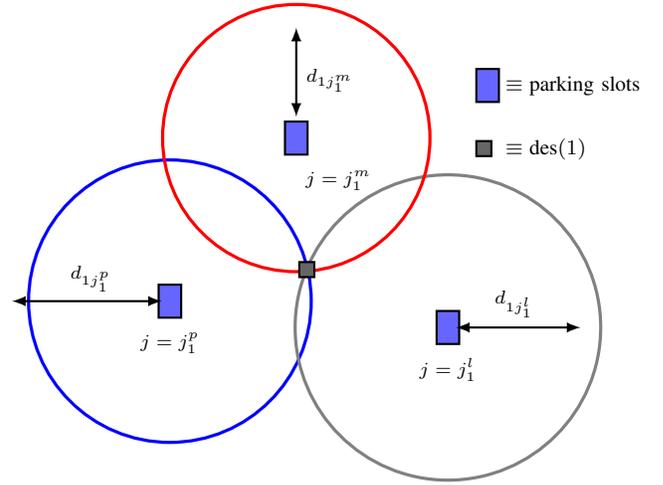
\begin{figure}[t]
\centering
\begin{tikzpicture}[scale=1.4]
\tikzstyle{ParkSlot}=[rectangle,rotate=90,draw=black,fill=blue!60!white,thick,inner sep=0pt,minimum size=3.0mm]
\tikzstyle{DesNode}=[rectangle,rotate=90,draw=black!100!white,fill=black!60!white,thick,inner sep=0pt,minimum size=2.0mm]
\tikzstyle{MyMatrix1} = [matrix, draw=blue, thick, fill=blue!20, text centered,font=\footnotesize,minimum height=2em]

\draw[color=blue,very thick] (0,0.75) circle (0.074\textwidth);
\draw[color=gray,very thick] (2.64,0.5) circle (0.08\textwidth);
\draw[color=red,very thick] (1.2, 2.30) circle (0.07\textwidth);

%
%

\node at (0,0.75) [ParkSlot] {\addb{.....}};
\node at (0,0.35) {\scriptsize{$j=j_1^{p}$}};

\node at (2.64,0.5) [ParkSlot] {\addb{.....}};
\node at (2.64,0.1) {\scriptsize{$j=j_1^{l}$}};

\node at (1.2,2.3) [ParkSlot] {\addb{.....}};
\node at (1.6,1.9) {\scriptsize{$j=j_1^{m}$}};

\node at (1.3,1.05) [DesNode] {};

\node at (2.9,2.2) [DesNode,anchor=north] {};
\node at (3.1,2.2) [anchor=west]{\footnotesize{$\equiv\mbox{des}(1)$}};

\node at (2.9,2.8)  [ParkSlot,anchor=north] {\addb{.....}};
\node at (3.1,2.8)[anchor=west] {\footnotesize{$\equiv\mbox{parking slots}$}};

\draw [latex-latex,black,thick,shorten >=0pt,shorten <=3pt] (0,0.75) -- (-1.5,0.75) node[pos=.5,above] {\scriptsize{$d_{1j_1^p}$}};
\draw [latex-latex,black,thick,shorten >=0pt,shorten <=3pt] (2.64,0.5) -- (3.9,0.5) node[pos=.5,above] {\scriptsize{$d_{1j_1^l}$}};
\draw [latex-latex,black,thick,shorten >=2pt,shorten <=8pt] (1.2, 2.30) -- (1.2, 3.40) node[pos=.5,right] {\scriptsize{$d_{1j_1^m}$}};

\end{tikzpicture}
\captionof{figure}{\label{fig:attack}\small{Given $(j_1^p,d_{1j_1^p})$, $(j_1^l,d_{1j_1^l})$, $(j_1^m,d_{ij_1^m})$ pairs known to the adversary (car~$n$), $j_1^p\neq j_1^l$, $j_1^p\neq j_1^m$, and $j_1^l\neq j_1^m$, discovering the location of $\mbox{des}(1)$.}}
\vspace{-5mm}
\end{figure}

Let us first fix the adversary to be car~$n$ and assume that car~$n$ wants to discover $\mbox{des}(1)$, i.e., the destination of car~$1$. Now \emph{suppose} car~$n$ knows the set $\mathcal{C}_1=\{(j_1^k,d_{1j_1^k})\}_{k=1,2,\ldots,K}$ of data associated with car~$1$, where $K$ is the total iterations of \textrm{DCP} algorithm. Provided there exists at least three distinct $j_1^k$s, car~$n$ can simply locate $\mbox{des}(1)$ as illustrated in Fig.~\ref{fig:attack}. Even \emph{if} car~$n$ knows only the set $\mathcal{D}_1=\{d_{1j_1^k}\}_{k=1,2,\ldots,K}$ of data associated with car~$1$, it turns out that car~$n$ can locate $\mbox{des}(1)$ exhaustively. In particular, in every iteration $k$, car~$n$ draws $M-1$ circles with radius $d_{1j_1^k}$ centered at parking slots $\mathcal{M}\setminus \{j_n^k\}$. Let $\mathcal{S}^k$ denote the aforementioned set of circles. Provided there are at least three distinct $j_1^k$s, which correspond to some iteration indexes $l,m$, and $p$, one can see that there exists at least one point at which a circle in $\mathcal{S}^l$, a circle in $\mathcal{S}^m$, and a circle in $\mathcal{S}^p$ intersect. If this point is unique, then it corresponds to $\mbox{des}(i)$.~\footnote{When the parking slots are not arbitrarily located and there are symmetric properties, then there can be more than one intersection point, which in turn will create uncertainties in correctly locating $\mbox{des}(i)$.} The discussion above indicates that \emph{if} the adversary (car~$n$) knows $\mathcal{C}_1$ or even $\mathcal{D}_1$, under mild conditions, it can locate $\mbox{des}(i)$. In the sequel, we show how \textrm{DCP} precludes such situations. In particular, we show how $\{d_{1j_1^k}\}_{k=1,2,\ldots,K}$ is kept hidden from the adversary car~$n$.

%
%
%

%
%

\subsection{How to preserve privacy} \label{subsec:Retain_Privacy}

Note that the only means by which car~$n$ gets access to some functions of $\{d_{1j_1^k}\}_{k=1,2,\ldots,K}$ is via $\{\lambda^{(k)}_n\}_{k=1,2,\ldots,K}$, see step~5 of \textrm{DCP} algorithm. In other words, the involvement of car~$n$ during the \textrm{DCP} algorithm is restricted so that, in every iteration~$k$, it has access to only some \emph{interface} variables $\lambda^{(k)}_n$ and ${\boldsymbol \mu}^{(k)}$.~\footnote{The knowledge of ${\boldsymbol \mu}^{(k)}$ is irrelevant here because it does not carry any information of $d_{1j_1^k}$.} This restriction is indeed achieved by the decomposition structure of problem~\eqref{eq:primal}. Moreover, we consider the situation that CC uses the step size~$\alpha_k$ of \textrm{DCP}~as
\be\label{eq:newStep}
\alpha_k = \frac{\alpha}{k} \ ,
\ee
where $\alpha$ is arbitrarily chosen on $[\alpha^{\textrm{min}},\alpha^{\textrm{max}}]$, $\alpha^{\textrm{min}}$ and $\alpha^{\textrm{max}}$ are positive numbers known only to CC such that $\alpha^{\textrm{min}}<\alpha^{\textrm{max}}$. Note that the above choice of $\alpha_k$ still preserves the convergence properties established in \emph{Proposition~1} (see Section~\ref{sec:convergence}) for the following reasons [compare with \eqref{eq:c11}]:
\begin{enumerate}
\item $\sum_{k=1}^{\infty} \alpha_k \rightarrow \infty$. This result is shown by noting that $\sum_{k=1}^{\infty}\alpha^{\textrm{min}}/k\rightarrow \infty$ and $0<\alpha^{\textrm{min}}/k \leq \alpha_k$ for all $k$.
\item $\sum_{k=1}^{\infty} \alpha^2_k$ converges to a point on $[(\alpha^{\textrm{min}})^2 \pi/6,(\alpha^{\textrm{max}})^2\pi/6]$. This result is achieved by noting that $\sum_{k=1}^{\infty}(\alpha^{\textrm{min}})^2/k^2\rightarrow(\alpha^{\textrm{min}})^2\pi/6$, $\sum_{k=1}^{\infty}(\alpha^{\textrm{max}})^2/k^2\rightarrow(\alpha^{\textrm{max}})^2\pi/6$, and $(\alpha^{\textrm{min}})^2/k^2\leq \alpha^2_k \leq (\alpha^{\textrm{max}})^2/k^2$ for all $k$.
\end{enumerate}
The arbitrary step size above (i.e., \ref{eq:newStep}) essentially introduces more protection to the problem data $\{d_{1j_1^k}\}_{k=1,2,\ldots,K}$.

Now we pose the following question: Does the proposed \textrm{DCP} allow car~$n$ alone to make records of $\{d_{1j_1^k}\}_{k=1,2,\ldots,K}$, so that it can locate $\mbox{des}(1)$ as discussed in Section~\ref{subsec:tent_attack}? It turns out that even though, car~$n$ can document the connections among the unknown parameters including $\{d_{1j_1^k}\}_{k=1,2,\ldots,K}$, among others, it can only come up with an \emph{under determined set of nonlinear equations}. Therefore $\{d_{1j_1^k}\}_{k=1,2,\ldots,K}$ cannot be computed as we will see next.

\begin{figure}[t]
\centering
\tikzset{hatchspread=3pt}
\tikzstyle{crosspoint}=[circle,draw=black,fill=blue, inner sep=0pt,minimum size=0.8mm]
\tikzstyle{crosspoint2}=[circle,draw=red,fill=red, inner sep=0pt,minimum size=0.8mm]
\tikzstyle{crosspoint3}=[circle,draw=white,fill=white, inner sep=0pt,minimum size=0.8mm]
\begin{tikzpicture}[scale=5]
\draw[-latex] (0,0) -- (1.5,0) node[anchor=north] {\scriptsize{$\lambda_1$}};
\draw	(0,0) node[anchor=north] {\scriptsize{$0$}}
		(1,0) node[anchor=north] {\scriptsize{$1$}}
        (0,1) node[anchor=east] {\scriptsize{$1$}};

\draw[dotted, draw=black,fill=red,opacity=.05,thick,dashed]
 (1,0) --
 (1.48,0.48) --
 (1.48,1.48) --
 (0.48,1.48) --
 (0,1) --
 (1,0) -- cycle ;
 \node at (0.85,1.3) [crosspoint3] (Xpoint1) {$A_1$};

\draw[very thick] (1,0) -- (0,1);
\draw[thin,opacity=.2] (0,1) -- (0.5,1.5);
\draw[thin,opacity=.2] (1,0) -- (1.5,0.5);
\draw[thick, dashed,opacity=.2] (0.25,0.75) -- (0.25,1.5);
\draw[-latex, dashed] (0.25,0.75) -- (0.37,1.45);
\draw[pattern=north west lines, pattern color=blue,opacity=.2]
 (0,1) --
 (0.48,1.48) --
 (0,1.48) --
 (0,1) --cycle ;
 \node at (0.37,1.45) [crosspoint2] (Xpoint1) {.};
\node at (0.12,1.3) [crosspoint3] (Xpoint1) {$A_2$};
\node at (0.33,1.45) (Xpoint1) {\scriptsize{${\vec a}$}};

\draw[pattern=north west lines, pattern color=blue,opacity=.2]
 (1,0) --
 (1.48,0) --
 (1.48,0.48) --
 (1,0) --cycle ;
\node at (1.3,0.12) [crosspoint3] (Xpoint1) {$A_3$};

\draw[very thick, dashed, draw= red,opacity=.3 ] (0.6,0.75) -- (1.35,1.5);
\draw[very thick, dashed, draw= red,opacity=.3 ] (1.1, 0.575) -- (1.5,0.975);
\node at (0.6,0.70) (Xpoint1) {\scriptsize{${\vec b}_1$}};
  \node at (1.1, 0.525) (Xpoint1) {\scriptsize{${\vec b}_2$}};
  \node at (0.6,0.75) [crosspoint] (Xpoint1) {.};
  \node at (1.1, 0.575) [crosspoint] (Xpoint1) {.};

\node at (0.25,0.75) [crosspoint] (Xpoint1) {.};
\node at (0.15,0.70) (Xpoint1) {\scriptsize{$\big(\lambda_1^{(1)},\lambda_2^{(1)}\big)$}};
\draw[-latex, very thick] (0.25,0.75) -- (0.85,1.0);

\node at (0.85,1.0) [crosspoint2] (Xpoint1) {.};
\node at (0.90,1.1) (Xpoint1) {\scriptsize{${\vec a}_1=\big(\lambda_1^{(1)}+\alpha_1d_{1j_1^1},\lambda_2^{(1)}+\alpha_1d_{2j_2^1}\big)$}};
\draw[thin,opacity=.2] (0.25,0.4) -- (0.6,0.75); 

\node at (0.425,0.575) [crosspoint] (Xpoint1) {.};
\node at (0.265,0.565) (Xpoint1) {\scriptsize{$\big(\lambda_1^{(2)},\lambda_2^{(2)}\big)$}};
\draw[-latex, very thick] (0.425,0.575) -- (1.2,0.675);
\node at (1.2,0.675) [crosspoint2] (Xpoint1) {.};
\node at (1.215,0.725) (Xpoint1) {\scriptsize{${\vec a}_2=\big(\lambda_1^{(2)}+\alpha_2d_{1j_1^2},\lambda_2^{(2)}+\alpha_2d_{2j_2^2}\big)$}};
\draw[thin,opacity=.2] (0.70,0.175) -- (1.1,0.575); 

\node at (0.7625,0.2375) [crosspoint] (Xpoint1) {.};
\node at (0.595,0.2275) (Xpoint1) {\scriptsize{$\big(\lambda_1^{(2)},\lambda_2^{(2)}\big)$}};

\draw[thin,dashed,opacity=.2] (0.425,0.575) -- (1.24,0.575);
\draw[thin,dashed,opacity=.2] (0.85,0.55) -- (0.85,1.05);
\draw[thin,dashed,opacity=.2] (0.25,0.75) -- (0.88,0.75);
\node at (0.6375,0.50) (Xpoint1) {$\underbrace{\addw{------}}_{\beta_1}$};

\draw[thin,dashed,opacity=.2] (0.7625,0.2375) -- (1.25,0.2375);
\draw[thin,dashed,opacity=.2] (1.2,0.220) -- (1.2,0.710);
\node at (0.98125,0.1625) (Xpoint1) {$\underbrace{\addw{------.}}_{\beta_2}$};
		

\draw[-latex] (0,0) -- (0,1.55) node[anchor=east] {\scriptsize{$\lambda_2$}};


\end{tikzpicture}
\captionof{figure}{\label{fig:lambda_evoltion}\small{The evolution of ${\boldsymbol\lambda}^{(k)}=(\lambda^{(k)}_1,\lambda^{(k)}_2)$ in a case of $2$-cars.}}
\vspace{-5mm}
\end{figure}
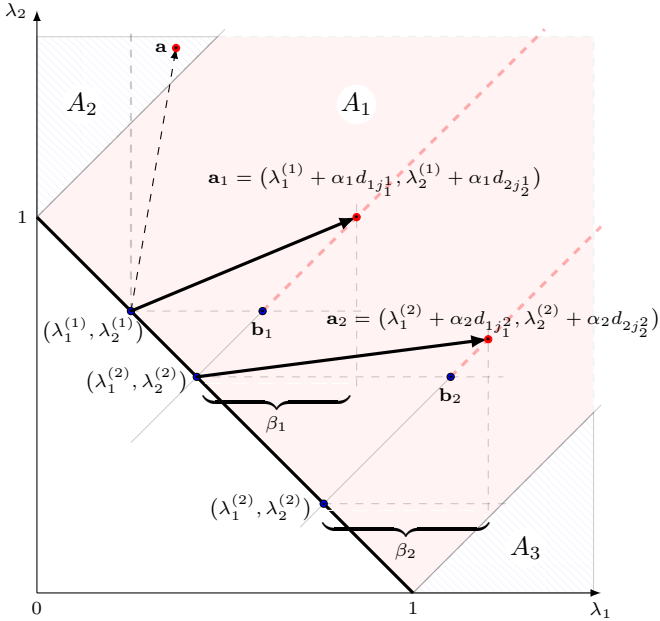

We consider only the case with $n=2$ and suppose car~$2$ is the adversary that wants to discover the destination of car~$1$, i.e., $\mbox{des}(1)$. The discussion can be generalized to scenarios with $n>2$, in a straightforward manner.

\begin{table*}[t]
\centering
\begin{tabular}{ c|r|c|c|c  }
\hline
iteration ($k$) & relations of unknown parameters \addw{abc} & unknowns & no. of unknowns ($U_k$) & no. of equations ($E_k$) \\ \hline
{$1$} & $\lambda^{(1)}_1+\lambda^{(1)}_2=1 \quad\rightarrow (1.1)$  & $\lambda^{(1)}_1$ & $1$ & $1$ \\ \hline
\multirow{4}{*}{$2$} & $(1.1)$ & $\lambda^{(1)}_1,\lambda^{(2)}_1$ & \multirow{4}{*}{5} & \multirow{4}{*}{4} \\
 & $\lambda^{(2)}_1+\lambda^{(2)}_2=1\quad\rightarrow (2.1)$ & $\beta_1$ & & \\
 & $\lambda^{(2)}_1+ \beta_1=\lambda^{(1)}_1+ \alpha_1 d_{1j_1^1}\quad\rightarrow (2.2)$ & $\alpha_1$ &  &  \\
 & $\lambda^{(2)}_2+ \beta_1=\lambda^{(1)}_2+ \alpha_1 d_{2j_2^1}\quad\rightarrow (2.3)$ & $d_{1j_1^1}$ &  &  \\ \hline
\multirow{4}{*}{$3$} & $(1.1),(2.1),(2.2),(2.3)$ & $\lambda^{(1)}_1,\lambda^{(2)}_1,\lambda^{(3)}_1$ & \multirow{4}{*}{9} & \multirow{4}{*}{7} \\
 & $\lambda^{(3)}_1+\lambda^{(3)}_2=1\quad\rightarrow (3.1)$ & $\beta_1,\beta_2$ & & \\
 & $\lambda^{(3)}_1+ \beta_2=\lambda^{(2)}_1+ \alpha_2 d_{1j_1^2}\quad\rightarrow (3.2)$ & $\alpha_1,\alpha_2$ &  &  \\
 & $\lambda^{(3)}_2+ \beta_2=\lambda^{(2)}_2+ \alpha_2 d_{2j_2^2}\quad\rightarrow (3.3)$ & $d_{1j_1^1},d_{1j_1^2}$ &  &  \\ \hline
{$\vdots$} & $\vdots$ & $\vdots$ & $\vdots$ & $\vdots$ \\ \hline
\end{tabular}
\caption{\small{Relations of unknown parameters as seen by the adversary car~$2$.}}
\label{tab:Car2Data}
\vspace{-5mm}
\end{table*}

First, note that CC performs the projection of ${\boldsymbol \lambda}^{(k)}-\alpha_k{\vec u}^{(k)}$ onto the probability simplex to yield ${\boldsymbol \lambda}^{(k+1)}$ [see step~4-c of \textrm{DCP} algorithm]. In the considered $2$-car case, the probability simplex is the line segment form $(0,1)$ to $(1,0)$, see Fig.~\ref{fig:lambda_evoltion}. Once car~$2$ is given $\lambda_2^{(1)}$, it can locate $(\lambda_1^{(1)},\lambda_2^{(1)})$, because $\lambda_1^{(1)}=1-\lambda_2^{(1)}$. Yet car~$2$ cannot locate ${\vec a}_1={\boldsymbol \lambda}^{(1)}-\alpha_1{\vec u}^{(1)}$. After receiving $\lambda_2^{(2)}$, car~$2$ can locate $(\lambda_1^{(2)},\lambda_2^{(2)})$. It can also locate ${\vec a}_1$ up to the \emph{ray} originating at ${\vec b}_1$, see Fig.~\ref{fig:lambda_evoltion}. However, car~$2$ cannot exactly locate ${\vec a}_1$. The algorithm continues in a similar manner. For example, the evolution of ${\boldsymbol{\lambda}}^{(k)}=(\lambda^{(k)}_1,\lambda^{(k)}_2)$ is illustrated in Fig.~\ref{fig:lambda_evoltion} for $k=1$, $2$, and $3$. With this knowledge of ${\boldsymbol\lambda}^{(k)}$ evolution, car~$2$ can write a set of equations in every iteration~$k$ as given in Table~\ref{tab:Car2Data}. Note that the set of equations for $k>1$ are nonlinear, because there are products of unknowns, e.g., $(2.2)$, $(2.3)$, $(3.2)$, and $(3.3)$.

When documenting the relations of unknowns in Table~\ref{tab:Car2Data}, we assume $\{{\vec a}_k\}_{k=1,2,3\ldots}$ lies only in the shaded area~$A_1$. In contrast, suppose ${\vec a}_k$ lies either in the hatched area $A_2$ or $A_3$ for some iterations. One such point is depicted in Fig.~\ref{fig:lambda_evoltion}, where ${\vec a}_1={\vec a}$, see the hatched area $A_2$. In this case, $(\lambda_1^{(2)},\lambda_2^{(2)})$ will be $(0,1)$. Consequently, car~$2$ can locate ${\vec a}_1$ only up to a \emph{cone} instead of a \emph{ray}, which in turn accounts for more uncertainties in determining ${\vec a_1}$. Such situations can only increase the difference between the number of unknowns ($U_k$) and the number of equations ($E_k$). For example, in this case, we will have $U_2-E_2>1$, instead of $U_2-E_2=1$ as in Table~\ref{tab:Car2Data}.

Thus, we can conclude that, the total unknowns are always greater than the total equations. This results an under determined set of nonlinear equations, see Table~\ref{tab:Car2Data}. Therefore, car~$2$ cannot make records of unknowns $\{d_{1j_1^k}\}_{k=1,2,\ldots,K}$ and consequently it cannot discover $\mbox{des}(1)$ as discussed in Section~\ref{subsec:tent_attack}.

\section{Numerical results } \label{sec:results}
In this section we present the numerical evaluation of our proposed algorithm DCP. We compare the \textrm{DCP} algorithm to the following benchmarks:
\begin{enumerate}
\item[(a)] Greedy parking policy: In this case, each car selects the closest parking slot to its destination.
\item[(b)] Optimal parking policy: A solution of the optimization problem~\eqref{eq:epigraph} is found by using the general solver, the IBM CPLEX optimizer~\cite{CPLEX}.
\end{enumerate}

In each time slot, the proposed algorithm is carried out for $K$ subgradient iterations. In fact, $K$ is used to define the stopping criterion at step~5 of \textrm{DCP} algorithm, see Fig.~\ref{fig:DAA}. In addition, the greedy policy and the optimal policy have also been performed at every time slot. We average the results over $T$ time slots to demonstrate the average performances of the \textrm{DCP} algorithm. Specifically, at the beginning of every time slot, the total number of cars~$N$ and the total number of free parking slots~$M$ are considered to be fixed and the distances $\{d_{ij}\}_{i\in\mathcal{N},j\in\mathcal{M}}$ are considered uniformly distributed on $[0,1000]$. The parking distances, $\{d_{ij}\}_{i\in\mathcal{N},j\in\mathcal{M}}$ are changed from slot to~slot.

To simplify the presentation, we denote by $p^{\textrm{cur}}(t,k)$ the best objective value achieved at time slot $t$ after $k$ subgradient iterations [compare to $p^{\textrm{cur}}(k)$ in step~4-c,d of the \textrm{DCP} algorithm]. In~particular,
\be
p^{\textrm{cur}}(t,k) = \arg\min_{l=1,\ldots,k}p(t,l) \ ,
\ee
where $p(t,l)$ is the objective value at time slot $t$ and at subgradient iteration $l$. Note that $p^{\textrm{cur}}(t,k)$ is similar to $p^{\textrm{cur}}(k)$ of \textrm{DCP} algorithm with an additional index $t$ to indicate the time slot. Moreover, we denote by ${\vec X}^{\textrm{cur}}(t,k)$ the \emph{best} feasible or infeasible solution, which is identical to ${\vec X}^{\textrm{cur}}(k)$ of \textrm{DCP} algorithm with an additional index $t$ to indicate the time slot.

In all the considered simulations, we use $T=1000$. Moreover, $K$ is chosen to be $300$ or $500$. To simplify the presentation, we refer to problem setups with $N/M\leq0.5$ as \emph{lightly loaded} cases and refer to problem setups with $N/M \geq0.5$ as \emph{heavily loaded} cases. Moreover, the problem setups with $N/M\simeq0.5$ are referred to as \emph{moderately loaded} cases.

We first define a performance metric called \emph{the degree of feasibility} of ${\vec X}^{\textrm{cur}}(t,k)$. Note that ${\vec X}^{\textrm{cur}}(t,k)$ is, in fact, the assignment at the beginning of step~6 of \textrm{DCP} algorithm, which can be either feasible or infeasible. If ${\vec X}^{\textrm{cur}}(t,k)$ is feasible, we have $N^{\textrm{conflict}}=0$ or equivalently, $p^{\textrm{cur}}(t,K)<\infty$ [compare with step~6, 4-c, and 4-d]. On the other hand, if ${\vec X}^{\textrm{cur}}(t,k)$ is infeasible, we have $N^{\textrm{conflict}}>0$ or equivalently, $p^{\textrm{cur}}(t,K)=\infty$ [compare with step~6, 4-c, and 4-d]. This motivates to define the degree of feasibility~(\textsc{df}) of ${\vec X}^{\textrm{cur}}(t,k)$ as
\be
\textsc{df}_{K} = \frac{\sum_{t=1}^{T}I\big(p^{\textrm{cur}}(t,K)<\infty\big)}{T}\times 100 \% \ ,
\ee
where $I(E)$ is the indicator function of event $E$, i.e., $I(E)=1$ if $E$ is true or $I(E)=0$~otherwise.

\begin{figure}[t]
\centering
\subfigure[]
{\includegraphics[width=0.47\textwidth]{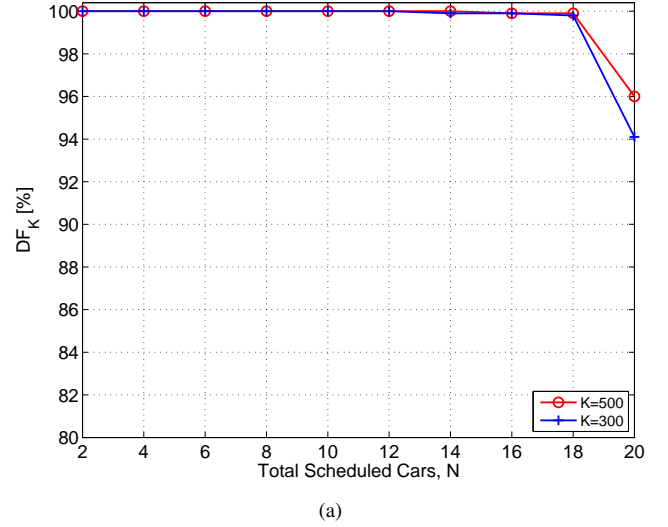}
\label{fig:DF_VS_N1}}
\goodgap
\subfigure[]
{\includegraphics[width=0.47\textwidth]{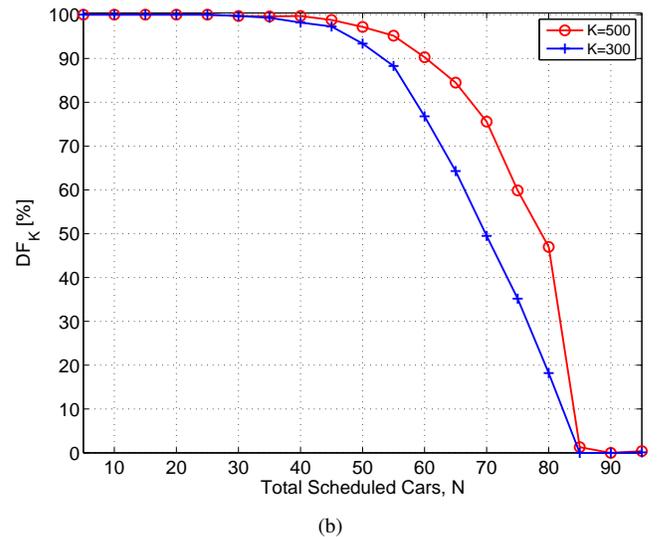}
\label{fig:DF_VS_N2}}\vspace{-2mm}
\caption{Degree of feasibility \textsc{df}$_K$ versus total cars $N$: (a) $M=20$ ; (b) $M=100$.}
\label{fig:DF_VS_N}
\vspace{-7mm}
\end{figure}

Fig.~\ref{fig:DF_VS_N} shows $\textsc{df}_K$ versus $N$ for fixed $M$. In particular a smaller dimensional problem with $M=20$~[Fig.~\ref{fig:DF_VS_N1}] and a larger dimensional problem with $M=100$~[Fig.~\ref{fig:DF_VS_N2}] is considered. Results show that when $N$ is significantly smaller than $M$, the degree of feasibility is almost $100\%$. Results further show that as $N$ becomes closer to $M$, the degree of feasibility starts deteriorating. Moreover, it becomes significantly law in the case of larger dimensional problem~[Fig.~\ref{fig:DF_VS_N2}] compared with the smaller dimensional problem~[Fig.~\ref{fig:DF_VS_N1}]. For example, when $N=M$, algorithm yields $\textsc{df}_K$ values in the range $94-96\%$ for the smaller dimensional problem setup. However, in the case of larger dimensional problem, when $N=M$, the degree of feasibility is almost zero. Not surprisingly, running \textrm{DCP} algorithm for larger number of subgradient iterations (e.g., $K=500$ ) yields better feasibility results compared with smaller number of subgradient iterations (e.g., $K=300$). However, the performance gap has been pronounced in the case of $M=100$~[Fig.~\ref{fig:DF_VS_N2}] compared to the case $M=20$.

\begin{figure}[t]
\centering
\subfigure[]
{\includegraphics[width=0.47\textwidth]{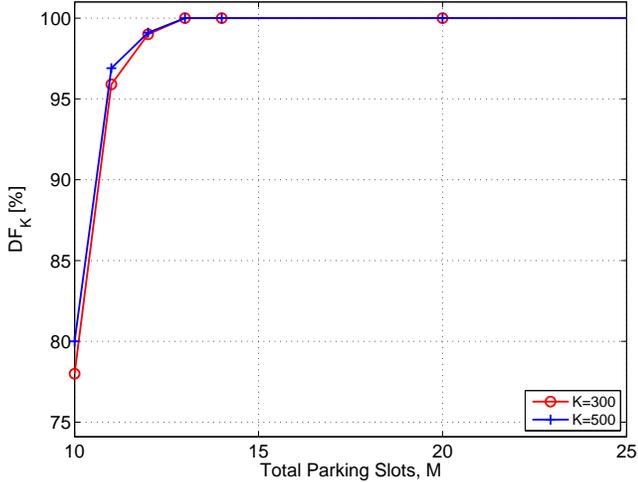}
\label{fig:DF_VS_M1}}
\goodgap
\subfigure[]
{\includegraphics[width=0.47\textwidth]{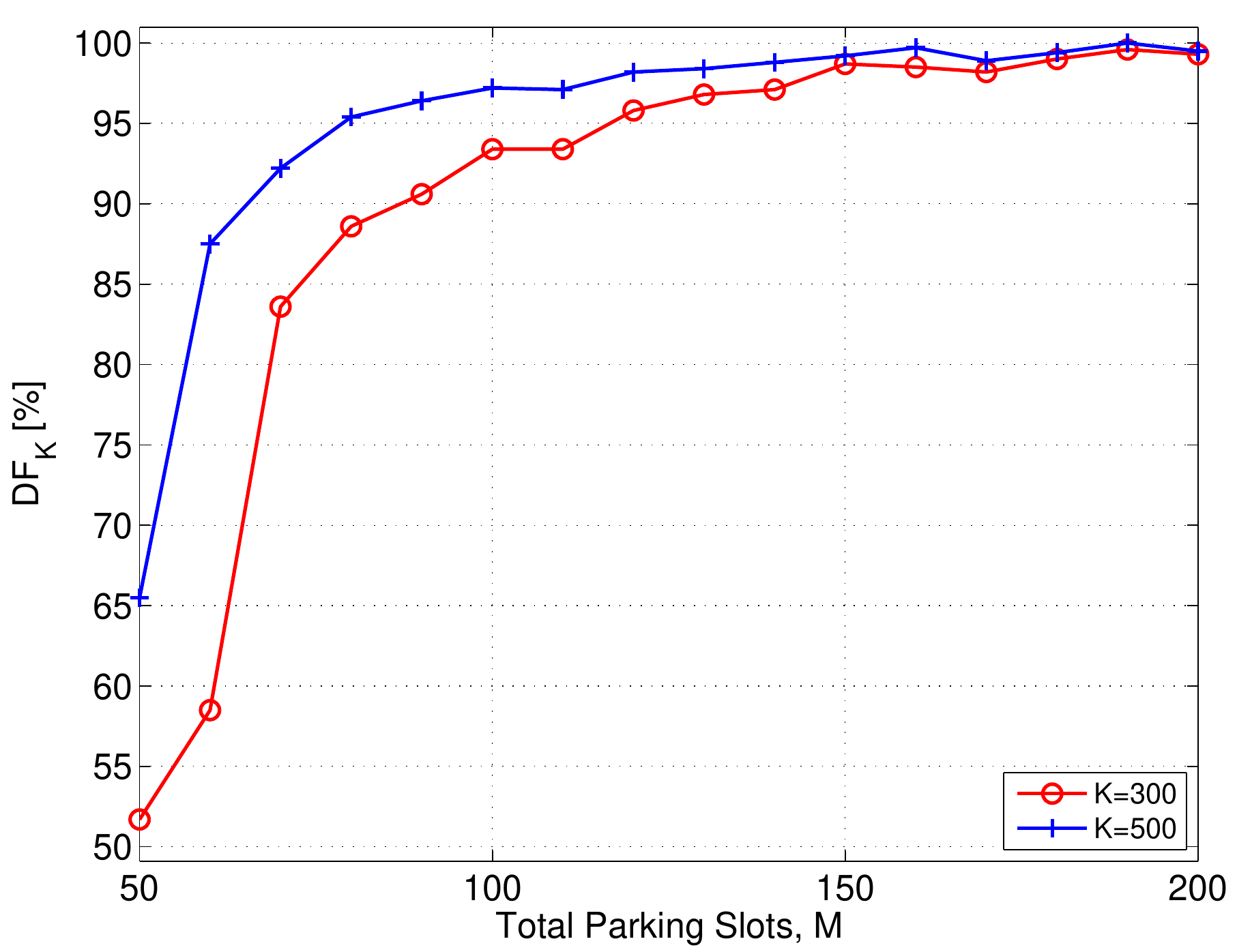}
\label{fig:DF_VS_M2}}\vspace{-2mm}
\caption{Degree of feasibility \textsc{df}$_K$ vs total parking slots $M$: (a) $N=10$ ; (b) $N=50$.}
\label{fig:DF_VS_M}
\vspace{-7mm}
\end{figure}

Fig.~\ref{fig:DF_VS_M} shows $\textsc{df}_K$ versus $M$ for fixed $N$. Again a smaller dimensional problem ($N=10$) and a larger dimensional problem with $N=50$ is considered, see Fig.~\ref{fig:DF_VS_M1}] and Fig.~\ref{fig:DF_VS_M2}, respectively. Results resemble the observations of Fig.~\ref{fig:DF_VS_N}, where a desirable feasibility is achieved when $M$ is significantly larger than $N$ and the performances are pronounced for smaller dimensional problems.

To see the average behavior of the \textrm{DCP} algorithm, now we consider the following performance~metric, which is a measure of the average objective value at subgradient iteration k:
\be\label{eq:p_avg}
p^{\textrm{ave}}(k)= \textstyle \frac{1}{T}\sum_{t=1}^{T}p^{\textrm{cur}}(t,k) \ , \qquad k=1,\ldots,K \ .
\ee

\begin{figure}[t]
\centering
\subfigure[]
{\includegraphics[width=0.47\textwidth]{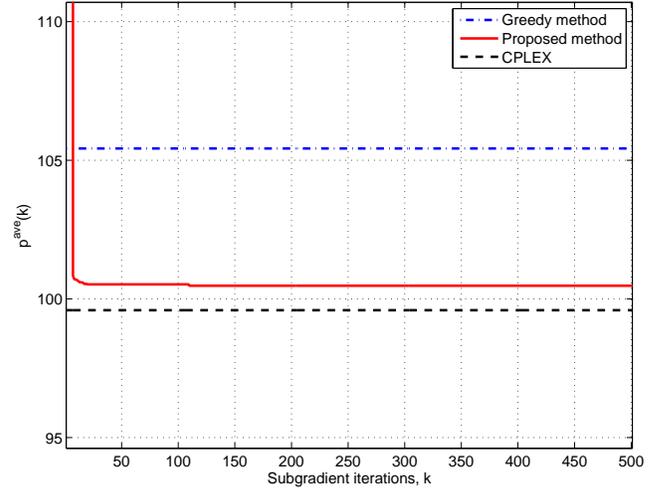}
\label{fig:AVG_obj_VS_k_small1}}
\goodgap
\subfigure[]
{\includegraphics[width=0.47\textwidth]{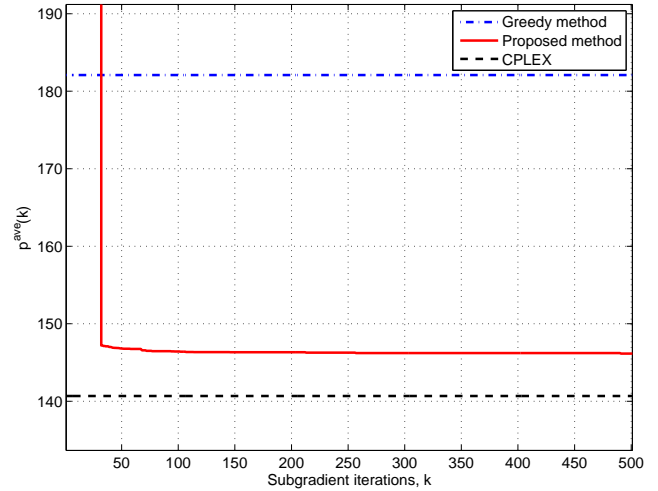}
\label{fig:AVG_obj_VS_k_small2}}\vspace{-2mm}
\caption{Average objective $p^{\textrm{ave}}(k)$ versus subgradient iterations $k$: (a) ${N}=4$ and ${M}=20$; (b)~${N}=10$ and ${M}=20$.}
\label{fig:AVG_obj_VS_k_small}
\vspace{-7mm}
\end{figure}
\begin{figure}[t]
\centering
\subfigure[]
{\includegraphics[width=0.47\textwidth]{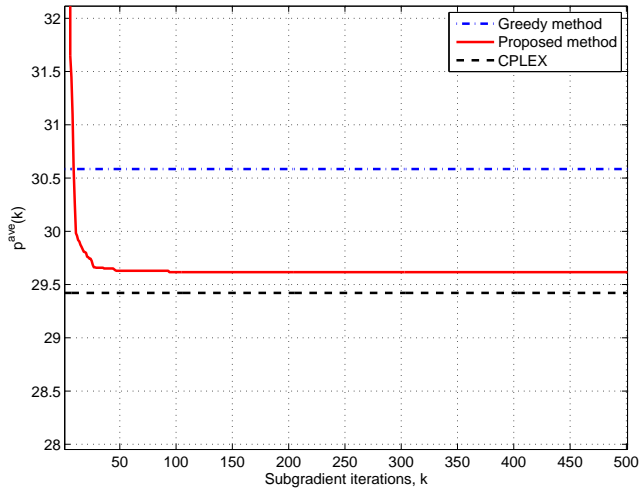}
\label{fig:AVG_obj_VS_k_large1}}
\goodgap
\subfigure[]
{\includegraphics[width=0.47\textwidth]{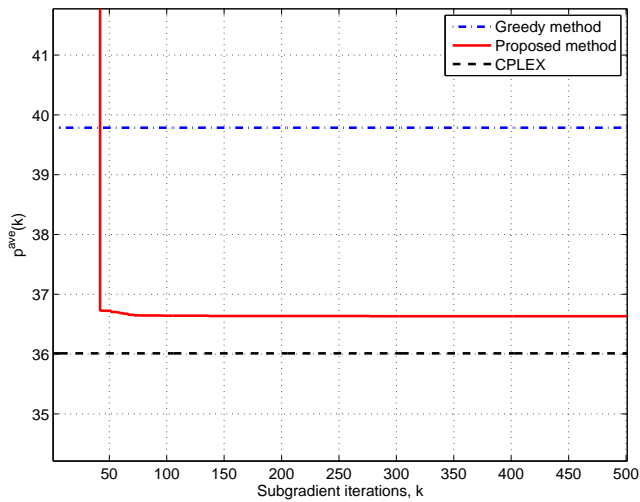}
\label{fig:AVG_obj_VS_k_large2}}\vspace{-2mm}
\caption{Average objective $p^{\textrm{ave}}(k)$ versus subgradient iterations $k$: (a) ${N}=10$ and ${M}=100$; (b)~${N}=20$ and ${M}=100$.}
\label{fig:AVG_obj_VS_k_large}
\vspace{-7mm}
\end{figure}

\begin{figure}[t]
\centering
\subfigure[]
{\includegraphics[width=0.47\textwidth]{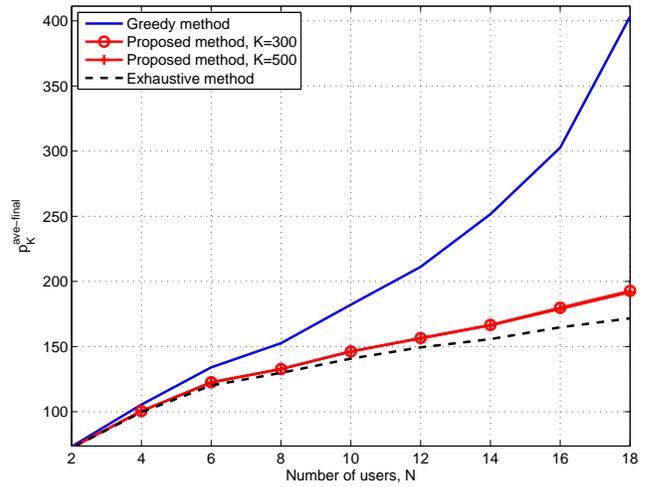}
\label{fig:AVG_Final_VS_N1}}
\goodgap
\subfigure[]
{\includegraphics[width=0.47\textwidth]{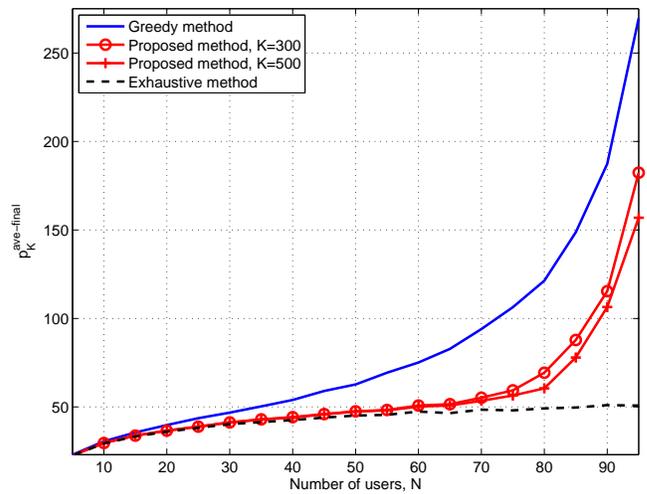}
\label{fig:AVG_Final_VS_N2}}\vspace{-2mm}
\caption{Average objective value after the termination of \textrm{DCP} $p^{\textrm{ave-final}}_K$ versus total cars~$N$: (a) $M = 20$; (b) $M = 100$.}
\label{fig:AVG_Final_VS_N}
\vspace{-7mm}
\end{figure}

Fig.~\ref{fig:AVG_obj_VS_k_small} shows $p^{\textrm{ave}}(k)$ versus subgradient iterations $k$ for cases $M=20$, $N=4$~[Fig.~\ref{fig:AVG_obj_VS_k_small1}] and $M=20$, $N=10$~[Fig.~\ref{fig:AVG_obj_VS_k_small2}]. Note that the vertical drops of the curves associated with our proposed method correspond to the subgradient iteration, before which a feasible assignment is found during any time slots $t=\{1,\ldots,T\}$. Results show that in the lightly loaded case that corresponds to a smaller $N$ (i.e., $N=4$), a feasible assignment is found much earlier than the moderately loaded case which corresponds to $N=10$. Specifically, when $N=4$, the \textrm{DCP} algorithm yields a feasible assignment at most after $k=6$ subgradient iterations, whereas when $N=10$, it yields a feasible assignment at most after $k=31$ subgradient iterations. This result is consistent with Fig.~\ref{fig:DF_VS_N} and Fig.~\ref{fig:DF_VS_M}, because for fixed $M$, the higher the $N$, the higher the number of time slots among $t=1,\ldots,T$ at which feasible solutions are achieved. For comparison, we also plot the average objective values obtained from the greedy method and the optimal CPLEX method. Not surprisingly, the optimal CPLEX method gives the best average objective, which is achieved at the expense of high computational complexity. However,
our proposed method trades off an increase in average objective value for a
low complexity in the algorithm, which gracefully scalable. Still the performance degradation of the proposed method is not critical. For example, the performance loss of \textrm{DCP} method compared with the optimal is $0.88\%$ in the case of $N=4$ and $3.89\%$ in the case of $N=10$. It is interesting to note that the proposed \textrm{DCP} algorithm outperforms the greedy method. For example, the performance degradation of the greedy method compared with the optimal is $5.86\%$ and $29.43\%$ in the cases of $N=4$ and $N=10$, respectively, which is significantly higher compared with our \textrm{DCP} method.

Fig.~\ref{fig:AVG_obj_VS_k_large} shows $p^{\textrm{ave}}(k)$ versus subgradient iterations $k$ for a larger dimensional problem setup. In particular, we consider the cases $M=100$, $N=10$~[Fig.~\ref{fig:AVG_obj_VS_k_large1}] and $M=100$, $N=20$~[Fig.~\ref{fig:AVG_obj_VS_k_large2}]. The behavior of the plots are similar to those in Fig.~\ref{fig:AVG_obj_VS_k_small}. In the case of $N=10$, the performance degradation of \textrm{DCP} method compared with the optimal is $0.66\%$ and that of the greedy method is $3.95\%$. In the case of $N=20$, the performance degradation of \textrm{DCP} method is $1.72\%$ and that of the greedy method is $10.48\%$. Results thus show that even in larger networks our proposed \textrm{DCP} method can outperform the greedy approach substantially.

As seen in Fig.~\ref{fig:DF_VS_N} and Fig.~\ref{fig:DF_VS_M}, as the number of cars $N$ becomes relatively closer to the total free parking slots $M$ (i.e., heavily loaded scenarios), infeasibility of ${\vec X}^{\textrm{cur}}(t,k)$ is usually inevitable. Therefore, as discussed in section~\ref{subsec:subroutine}, we have to rely on a subroutine to construct a feasible assignment, see step~6 in \textrm{DCP} algorithm. In the sequel, we show numerically the performance of proposed \textrm{DCP} algorithm in such heavily loaded cases.

Recall that, if $p^{\textrm{cur}}(t,K)=\infty$, then the corresponding assignment $X^{\textrm{cur}}(t,K)$ is infeasible. In such situations, our proposed algorithm invokes its subroutine (section~\ref{subsec:subroutine}) to construct a feasible assignment by using the best infeasible solution achieved so far, see step~6 of \textrm{DCP}. On the other hand, if $p^{\textrm{cur}}(t,K)<\infty$, then the corresponding assignment $X^{\textrm{cur}}(t,K)$ is already feasible. In either case, \textrm{DCP} algorithm returns a feasible point as given in step~6 of \textrm{DCP} algorithm. We denote by $X^{\textrm{final}}(t)$ this feasible assignment returned by our proposed \textrm{DCP} at time slot $t$ and by $p^{\textrm{final}}(t,K)$ the corresponding objective value. Finally, we denote by $p^{\textrm{ave-final}}_K$ the average objective value achieved after the termination of \textrm{DCP} algorithm. In particular, $p^{\textrm{ave-final}}_K$ is given by
\be
p^{\textrm{ave-final}}_K = \textstyle \frac{1}{T}\sum_{t=1}^{T} p^{\textrm{final}}(t,K) \ .
\ee
Note that when $X^{\textrm{cur}}(t,K)$ is feasible for all $t\in\{1,\ldots,T\}$, then $p^{\textrm{ave-final}}_K=p^{\textrm{ave}}(K)$ [compare with \eqref{eq:p_avg}].

Fig.~\ref{fig:AVG_Final_VS_N} shows average objective value $p^{\textrm{ave-final}}_K$ versus $N$ for the cases $M=20$~[Fig.~\ref{fig:AVG_Final_VS_N1}] and $M=100$~[Fig.~\ref{fig:AVG_Final_VS_N2}]. Results show that \textrm{DCP} algorithm always outperforms the greedy method. Using subgradient iterations $K=500$ accounts for an increase in the performance gain compared with $K=300$, though the gains are not substantial. When the setup is lightly loaded, proposed \textrm{DCP} performs very close to the optimal approach. For fixed $M$, increasing $N$ results increasing of the performance gap. For example, in the case of $M=20$ and $N=8$ with $K=300$ corresponds to a $3.94\%$ performance deviation of \textrm{DCP} compared to the optimal CPLEX and $N=18$ corresponds to a $23.14\%$ performance deviation, see Fig.~\ref{fig:AVG_Final_VS_N1}. Moreover, in the case of $M=100$ and $N=50$ with $K=300$ corresponds to a $5.40\%$ performance deviation of \textrm{DCP} compared to the optimal CPLEX and $N=95$ corresponds to a $259.20\%$ performance deviation, Fig.~\ref{fig:AVG_Final_VS_N2}. Such reductions in the performance gains are certainly expected, because there is a trade-off between the complexity of the algorithms and the performance loss. Results further show that the larger the parking slots $M$, the larger the performance deviation, especially in heavily loaded cases.


\begin{figure}[t]
\centering
\subfigure[]
{\includegraphics[width=0.47\textwidth]{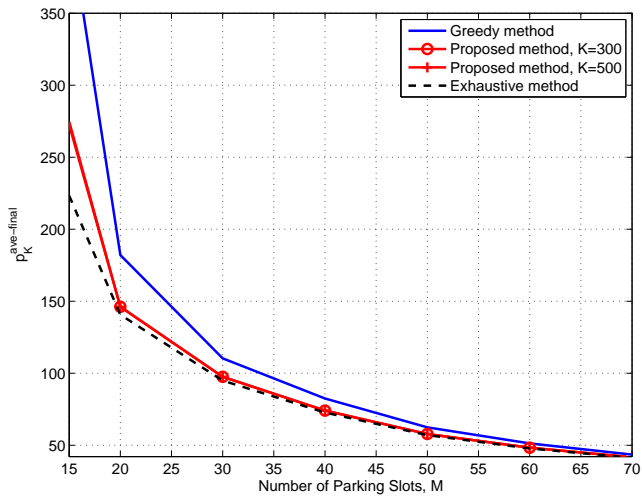}
\label{fig:AVG_Final_VS_M1}}
\goodgap
\subfigure[]
{\includegraphics[width=0.47\textwidth]{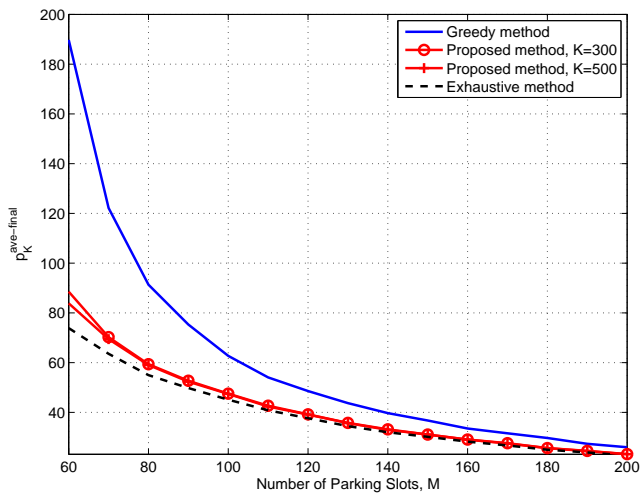}
\label{fig:AVG_Final_VS_M2}}\vspace{-2mm}
\caption{Average objective value after the termination of \textrm{DCP} $p^{\textrm{ave-final}}_K$ versus total parking slots~$M$: (a) $N = 10$; (b) $N = 50$.}
\label{fig:AVG_Final_VS_M}
\vspace{-7mm}
\end{figure}

Fig.~\ref{fig:AVG_Final_VS_M} shows average objective value $p^{\textrm{ave-final}}_K$ versus $M$ for the cases $N=10$~[Fig.~\ref{fig:AVG_Final_VS_M1}] and $N=50$~[Fig.~\ref{fig:AVG_Final_VS_M2}]. The observations are similar to those in Fig.~\ref{fig:AVG_Final_VS_N}. Results confirm that, our \textrm{DCP} algorithm performs very close to the optimal CPLEX method in lightly loaded cases, where $N/M<0.5$, see Fig.~\ref{fig:AVG_Final_VS_M1} curves with $M>20$ and Fig.~\ref{fig:AVG_Final_VS_M2} curves with $M>100$.  However, there is a noticeable performance degradation in heavily loaded cases, see Fig.~\ref{fig:AVG_Final_VS_M1} curves with $M<20$ and Fig.~\ref{fig:AVG_Final_VS_M2} curves with $M<100$. Moreover, the proposed method substantially outperforms the greedy method in all considered scenarios.

\begin{figure}[t]
\centering
\subfigure[]
{\includegraphics[width=0.47\textwidth]{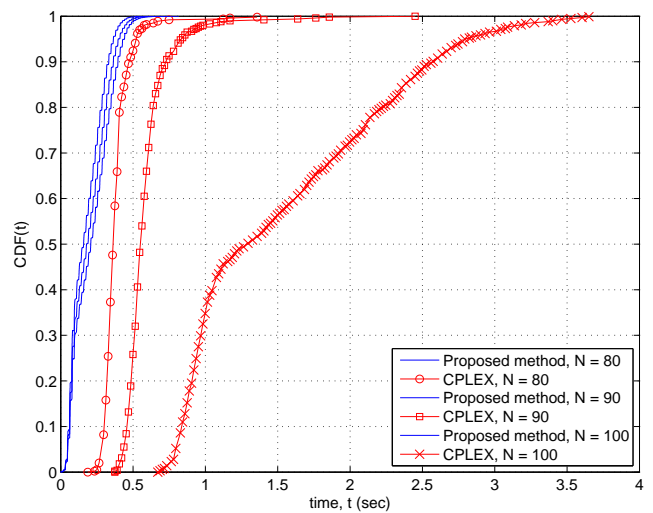}
\label{fig:CDF_plot1}}
\goodgap
\subfigure[]
{\includegraphics[width=0.47\textwidth]{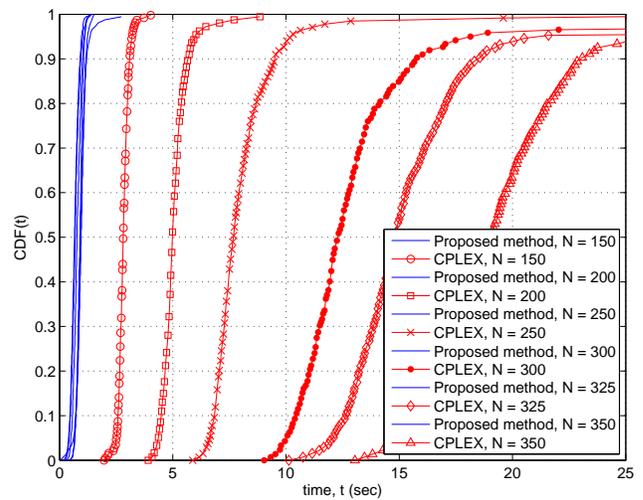}
\label{fig:CDF_plot2}}\vspace{-2mm}
\caption{CDF of time: (a) ${M}=100$; (b) ${M} = 500$.}
\label{fig:CDF_plot}
\vspace{-7mm}
\end{figure}

In order to provide a statistical description of the speed of
the proposed algorithm, we consider empirical the cumulative distribution function~(CDF) plots. Specifically, for each time slot $t\in\{1,\ldots,{T}\}$, we store the total CPU time required for \textrm{DCP} to find $X^{\textrm{final}}(t,K)$, where we use $K=300$. Similarly, the total CPU time required to find the \emph{optimal} value by using CPLEX is recorded. Figure~\ref{fig:CDF_plot} shows the empirical CDF plots of the time for $N=80,90,100$ with $M=100$~[Figure~\ref{fig:CDF_plot1}] and for $N=150,200,250,300,325,350$ with $M=500$~[Figure~\ref{fig:CDF_plot2}]. In the case of \textrm{DCP} algorithm, the effects of changing the problem size by increasing $N$ on the CDF plots are almost indistinguishable. However, in the case of optimal CPLEX method, there is a prominent increase in the time required to compute the optimal value. It should be emphasized that CPLEX finds the optimal assignment with a hight penalty on time to do it, where as \textrm{DCP} algorithm finds feasible assignment efficiently with a penalty on the optimality.

\begin{figure}[t]
\centering
\subfigure[]
{\includegraphics[width=0.47\textwidth]{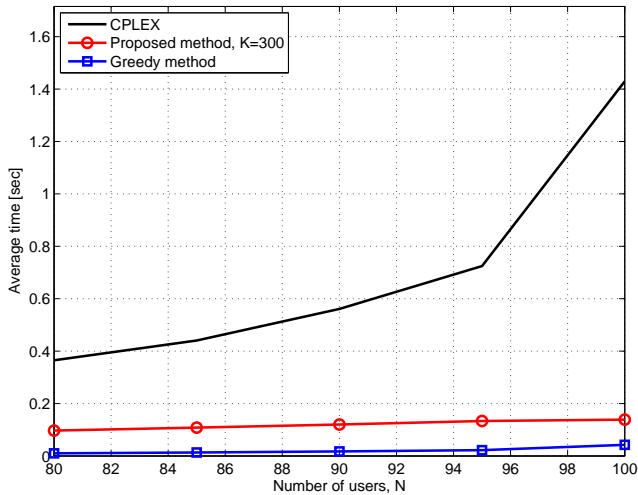}
\label{fig:AverageTime_plot1}}
\goodgap
\subfigure[]
{\includegraphics[width=0.47\textwidth]{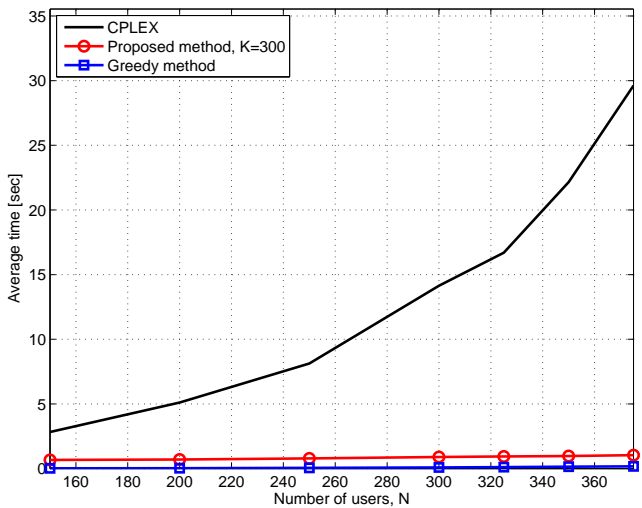}
\label{fig:AverageTime_plot2}}\vspace{-2mm}
\caption{Average CPU time: (a) ${M}=100$; (b) ${M} = 500$.}
\label{fig:AverageTime_plot}
\vspace{-7mm}
\end{figure}

Figure~\ref{fig:AverageTime_plot} depicts the average time required by DCP algorithm, the optimal method, and the greedy method versus $N$ for $M=100$~[Figure~\ref{fig:AverageTime_plot1}] and for $M=500$~[Figure~\ref{fig:AverageTime_plot2}]. Results show that the average time required by DCP to find possibly a suboptimal solution is not sensitive to the variation of $N$ and is in the range $0-1$ seconds. Moreover, they are comparable to the average time of simple greedy method. However, the average time required by the optimal method to find the optimal solution grows approximately exponentially with $N$. This is certainly expected because problem~(\ref{eq:epigraph}) is combinatorial, and therefore the worst-case complexity of the optimal method grows exponentially with the problem size~\cite[\S~1.4.2]{Boyd-Vandenberghe-04}. Thus, there is naturally a tradeoff between the optimality and the efficiency of the algorithms. The results together with those of Fig.~\ref{fig:AVG_Final_VS_N} and \ref{fig:AVG_Final_VS_M} suggest that our proposed DCP algorithm yields a good tradeoff between the optimality and the efficiency, especially in lightly loaded cases. These properties are favorable for practical implementation.

\section{Conclusions}\label{sec:conclusions}
In this paper, we considered the problem of car parking assignment. Unlike the existing greedy approaches, our problem formulation considered fairness among the scheduled cars in the sense that the global objective was to minimize
the maximum distance from the parking slots to the intended destinations
of the cars. A method based on Lagrange duality theory was proposed to address the nonconvex and combinatorial assignment problem. Even though we placed a stronger emphasis in the car parking slot assignment problem, our formulation and the corresponding algorithm generally applies in fair agent-target assignment problems in other application domains as well. We highlighted appealing privacy properties of the proposed algorithm. In particular, we showed that the proposed method is privacy preserving in the sense that any car involved in the algorithm will not be able to discover the destination of any other car during the algorithm iterations. Unlike the optimal exponentially complex approaches, our proposed method is scalable. Roughly speaking, numerical results showed that for all considered cases, where the number of free parking slots are equal or higher than twice the scheduled cars, our proposed algorithm's performance is similar to that of the optimal method. Moreover, in all considered cases, the proposed algorithm outperformed the simple greed approach. Therefore, the proposed algorithm yields a good trade-off between the implementation-level simplicity and the optimality.

\appendices
\setcounter{equation}{0}

\section{The Euclidean projection onto the probability simplex}\label{app:ProjSimplex}
\numberwithin{equation}{section}

In this appendix, we show how to project an arbitrary vector ${\vec x}\in\R^N$ onto the probability simplex~\eqref{eq:simplex}. This problem can formally be expressed as
\begin{IEEEeqnarray}{lcl}\label{eq:ProjSimplex}
\mbox{minimize} & \ \ & (1/2)||\boldsymbol \lambda-{\vec x}||^2_2 \IEEEyessubnumber\label{eq:ProjSimplex1}\\
\mbox{subject to} & \ \ &  {\vec 1}\tran \boldsymbol\lambda = 1 \IEEEyessubnumber\label{eq:ProjSimplex2}\\
& \ \  & \boldsymbol\lambda \succeq {\vec 0} \IEEEyessubnumber\label{eq:ProjSimplex5} \ ,
\end{IEEEeqnarray}
where the variable is $\boldsymbol\lambda\in\R^N$. Here ${\vec 1}$ denotes the vector with all $1$~s, ${\vec 0}$ denotes the vector with all $0$~s, and $\succeq$ denotes the component-wise inequality. The projection onto the probability simplex~\eqref{eq:simplex} is given by the solution $\bar{\boldsymbol\lambda}=(\bar\lambda_1,\ldots,\bar\lambda_N)$ of problem~\eqref{eq:ProjSimplex}. We apply duality theory~\cite[\S~5]{Boyd-Vandenberghe-04} to solve~\eqref{eq:ProjSimplex}.

To do this, we first form the partial Lagrangian by dualizing the constraint~\eqref{eq:ProjSimplex2}. Let $\nu$ denote the associated multiplier. Then the partial Lagrangian $L(\nu, \lambda)$ is given by
\begin{equation}
L(\nu, \boldsymbol\lambda) = (1/2)||\boldsymbol \lambda-{\vec x}||^2_2 + \nu({\vec 1}\tran\boldsymbol\lambda-1) \ .
\end{equation}
The dual function $h(\nu)$ is given by
\begin{subequations}\label{eq:dual_function20}
\begin{align}\label{eq:dual_function21}
   \hspace{-0mm}&  h(\nu) {=} \inf_{\substack{\lambda \succeq {\vec 0} }}L(\nu) {=} \inf_{\substack{\lambda \succeq {\vec 0} }}(1/2)||\boldsymbol \lambda{-}{\vec x}||^2_2 {+} \nu({\vec 1}\tran\boldsymbol\lambda{-}1)\\ \label{eq:dual_function_22}
  & =  \inf_{\substack{\lambda \succeq {\vec 0} }}(1/2)||\boldsymbol \lambda{-}({\vec x}-\nu{\vec 1})||^2_2 {+} \nu({\vec 1}\tran {\vec x}-1) {-} (1/2)N\nu^2  \\ \label{eq:dual_function_23}
  & = (1/2)\textstyle\sum_{i=1}^N \big(\min\{0,x_i{-}\nu\}\big)^2 {+} \nu({\vec 1}\tran {\vec x}{-}1) {-} (1/2)N\nu^2 ,
\end{align}
\end{subequations}
where the equality~\eqref{eq:dual_function_22} follows from straightforward mathematical manipulations and \eqref{eq:dual_function_23} follows from that $\lambda_i=0$ if $(x_i-\nu)\leq 0$ and $\lambda_i=(x_i-\nu)$ if $(x_i-\nu)> 0$. The dual optimal value $\nu^\star$ is given by $\nu^\star= \arg\max_{\nu}h(\nu)$~\cite[\S~5.2]{Boyd-Vandenberghe-04}

Note that the dual function $h(\nu)$ is a scalar valued function. Moreover, a subgradient $r(\nu)$ of $h(\nu)$ at $\nu$ can be analytically computed as
\be\label{eq:subgradient_projection}
r(\nu)= -\textstyle\sum_{i=1}^N \big(\min\{0,x_i-\nu\}\big) + {\vec 1}\tran {\vec x}-1 - N\nu \ .
\ee
Note that the dual function is always concave. Thus the sign of its subgradients changes as we pass through the maximum point (i.e., $\nu^\star$) of the dual function, which allows us to use a bisection search method to find $\nu^\star$ as follows.

\noindent\rule{0.49\textwidth}{0.3mm}
\\
\emph{Algorithm}: \ \textsc{\small{Computation of} \ $\nu^\star$}
\begin{enumerate}
\item Given the accuracy level $\epsilon>0$, $\nu^{\textrm{min}}$ and $\nu^{\textrm{max}}$ such that $r(\nu^{\textrm{min}})>0$ and $r(\nu^{\textrm{max}})<0$.
\item If $\nu^{\textrm{max}}-\nu^{\textrm{min}}<\epsilon$, return $\nu^\star=(\nu^{\textrm{max}}+\nu^{\textrm{min}})/2$ and STOP.
\item Set $\nu=(\nu^{\textrm{max}}+\nu^{\textrm{min}})/2$.
\item If $r(\nu)\geq0$, set $\nu^{\textrm{min}}=\nu$. Otherwise, set $\nu^{\textrm{max}}=\nu$. Go to step~2.
\end{enumerate}
\vspace{-4mm}
\rule{0.49\textwidth}{0.3mm}\vspace{-0mm}

Because strong duality holds for problem~\eqref{eq:ProjSimplex}~\cite[\S~5.2.3]{Boyd-Vandenberghe-04}, we can recover its solution $\bar{\boldsymbol\lambda}=(\bar\lambda_1,\ldots,\bar\lambda_N)$ as
\be
\bar\lambda_i = \max\{0, x_i-\nu^\star\} \ , \ i=1,\ldots,N \ .
\ee

\bibliographystyle{IEEEbib}
\bibliography{jour_short,conf_short,references,referencesElisabetta_04_22}

\begin{biography}[{\includegraphics[width=1in,height=1.25in,clip]{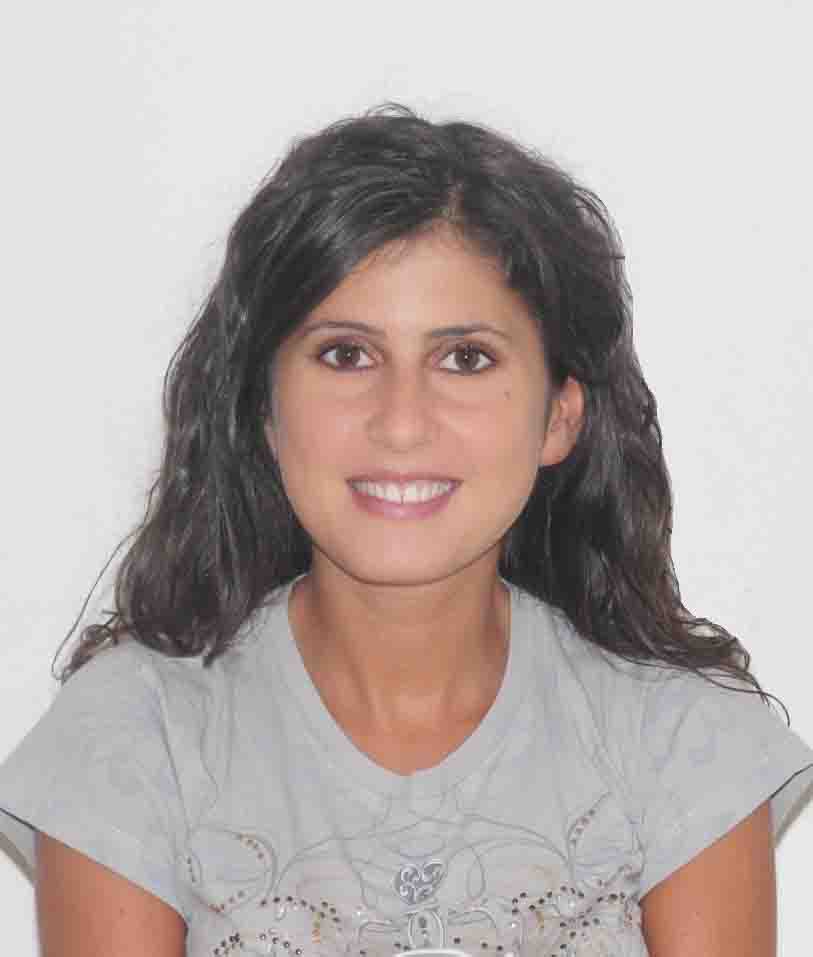}}]{Elisabetta Alfonsetti} received the Master degree in Computer Science from University of L'aquila, Italy, in 2012. Then she worked as research engineer in the Automatic Control Lab, Electrical Engineering Department and ACCESS Linnaeus Center, KTH Royal Institute of Technology, Stockholm, Sweden, focusing on optimization techniques and privacy preserving issue. She is currently with the TerraSwarm Lab, Electrical Engineering Department, UC Berkeley, California. Her research interests include application of Contract-based design paradigm for the design of complex systems.
\end{biography}
\begin{biography}[{\includegraphics[width=1in,height=1.25in,clip,keepaspectratio]{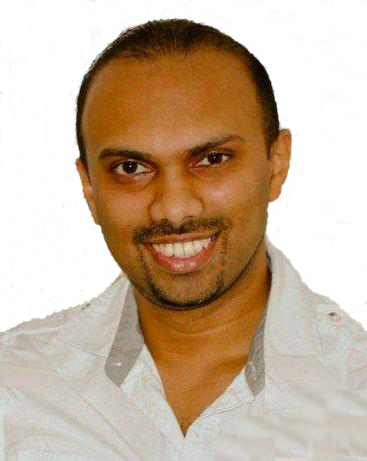}}]{Pradeep Chathuranga Weeraddana} (S'08,M'11)
received the M.Eng degree Telecommunication from School of Engineering and Technology, Asian Institute of Technology, Thailand in 2007 and the Ph.D. degree from University of Oulu, Finland, in 2011. He is currently working as postdoctoral researcher in Automatic Control Lab, Electrical Engineering Department and ACCESS Linnaeus Center, KTH Royal Institute of Technology, Stockholm, Sweden. His research interests include application of optimization techniques in various application domains, such as signal processing, wireless communications, smart grids, privacy, and security.
\end{biography}
\begin{biography}[{\includegraphics[width=1in,height=1.25in,clip,keepaspectratio]{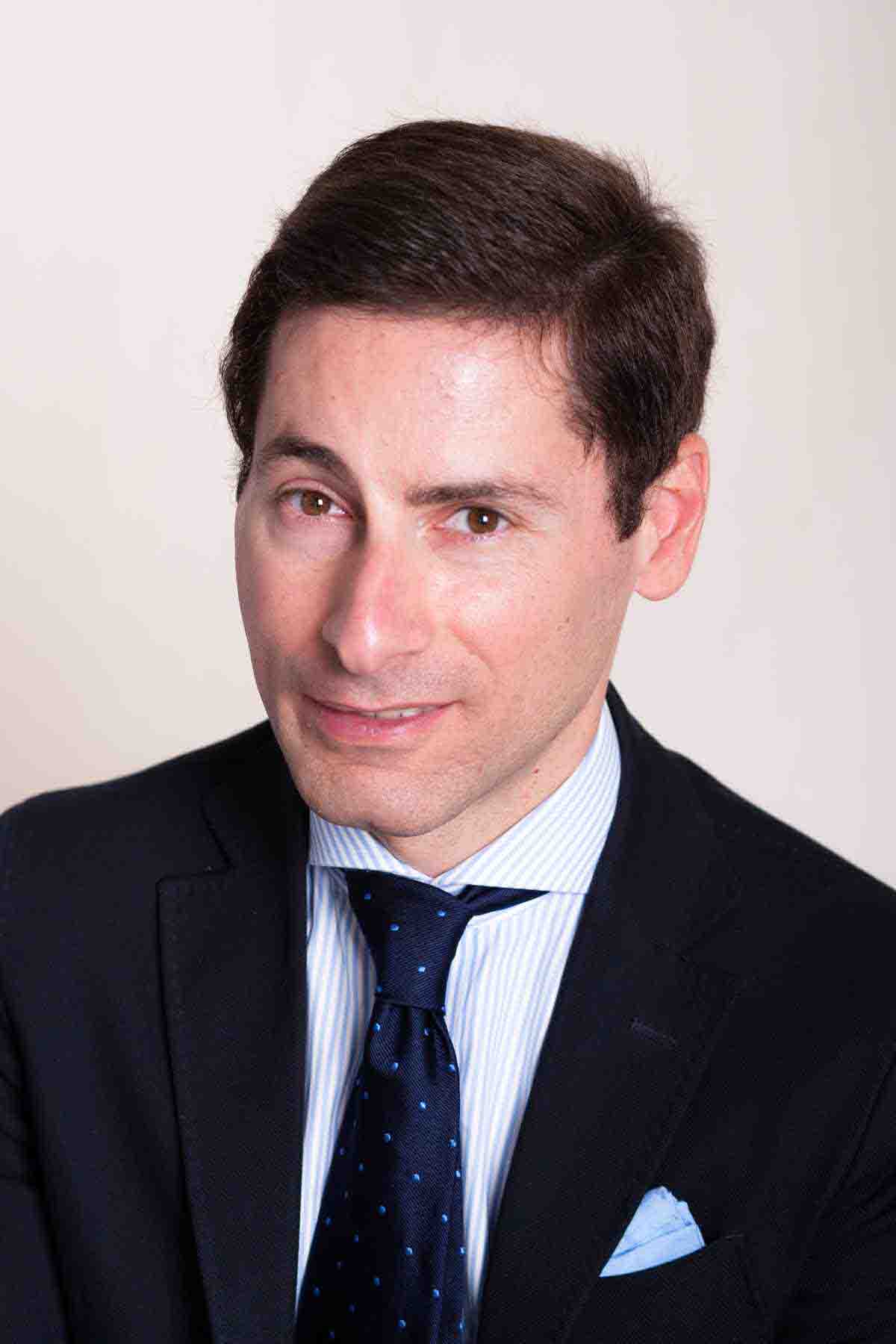}}]{Carlo Fischione} (M'05) is a tenured associate professor at KTH Royal Institute of Technology, Electrical Engineering, Stockholm, Sweden. His research interests include optimization and parallel computation with applications to wireless sensor networks, networked control systems, and wireless networks. He received the best paper award from the IEEE Transactions on Industrial Informatics of 2007. He is a member of the IEEE and SIAM.
\end{biography}

\end{document}

%% file: commands.tex
\newtheorem{example}{Example}

\newcommand{\be}{\begin{equation}}
\newcommand{\ee}{\end{equation}}
\newcommand{\ba}{\begin{array}}
\newcommand{\ea}{\end{array}}
\newcommand{\bea}{\begin{eqnarray}}
\newcommand{\eea}{\end{eqnarray}}

\newcommand{\tran}{^{\mbox{\scriptsize T}}}  

\newcommand{\ssum}[1]{\displaystyle\mathop{ \textstyle{\sum}}_{#1}}

\newcommand{\vbar}{\raisebox{.17ex}{\rule{.04em}{1.35ex}}}
\newcommand{\vbarind}{\raisebox{.01ex}{\rule{.04em}{1.1ex}}}
\newcommand{\R}{\ifmmode {\rm I}\hspace{-.2em}{\rm R} \else ${\rm I}\hspace{-.2em}{\rm R}$ \fi}
\newcommand{\D}{\ifmmode {\rm I}\hspace{-.2em}{\rm D} \else ${\rm I}\hspace{-.2em}{\rm D}$ \fi}
\newcommand{\T}{\ifmmode {\rm I}\hspace{-.2em}{\rm T} \else ${\rm I}\hspace{-.2em}{\rm T}$ \fi}
\newcommand{\N}{\ifmmode {\rm I}\hspace{-.2em}{\rm N} \else \mbox{${\rm I}\hspace{-.2em}{\rm N}$} \fi}
\newcommand{\B}{\ifmmode {\rm I}\hspace{-.2em}{\rm B} \else \mbox{${\rm I}\hspace{-.2em}{\rm B}$} \fi}
\newcommand{\Hil}{\ifmmode {\rm I}\hspace{-.2em}{\rm H} \else \mbox{${\rm I}\hspace{-.2em}{\rm H}$} \fi}
\newcommand{\C}{\ifmmode \hspace{.2em}\vbar\hspace{-.31em}{\rm C} \else \mbox{$\hspace{.2em}\vbar\hspace{-.31em}{\rm C}$} \fi}
\newcommand{\Cind}{\ifmmode \hspace{.2em}\vbarind\hspace{-.25em}{\rm C} \else \mbox{$\hspace{.2em}\vbarind\hspace{-.25em}{\rm C}$} \fi}
\newcommand{\Q}{\ifmmode \hspace{.2em}\vbar\hspace{-.31em}{\rm Q} \else \mbox{$\hspace{.2em}\vbar\hspace{-.31em}{\rm Q}$} \fi}
\newcommand{\Z}{\ifmmode {\rm Z}\hspace{-.28em}{\rm Z} \else ${\rm Z}\hspace{-.28em}{\rm Z}$ \fi}

\renewcommand{\vec}[1]{\bf{#1}}     

